\documentclass[12pt,authoryear]{elsarticle}
\usepackage[T1]{fontenc}
\usepackage[utf8]{inputenc}
\usepackage{lmodern}
\usepackage[a4paper,margin=2.5cm]{geometry}
\usepackage{amsmath,amssymb,amsthm}
\usepackage{graphicx}
\usepackage{booktabs,array,longtable}
\usepackage{subcaption}
\usepackage{placeins,float}
\usepackage{setspace}
\usepackage{lineno}
\usepackage{xurl}
\usepackage[hidelinks]{hyperref}
\biboptions{authoryear,round}
\graphicspath{{./}{Figures/}}
\allowdisplaybreaks
\newcommand{\horizontal}{horizontal-fairness}
\newcommand{\horizontalVertical}{horizontal-vertical-fairness}
\newcommand{\utilitarian}{utilitarian}
\newtheorem{proposition}{Proposition}
\newtheorem{lemma}{Lemma}

\newif\ifblindreview
\blindreviewfalse

\hypersetup{pdfauthor={},pdftitle={Restoring Power Fairly: Equitable Resource Allocation in Post-Disaster Grid Recovery}}

\begin{document}
\begin{frontmatter}
\title{Restoring Power Fairly: Equitable Resource Allocation in Post-Disaster Grid Recovery}
\ifblindreview
\else
\author[toronto]{Jie Mao}
\ead{jane.mao@mail.utoronto.ca}
\author[fau]{Saharnaz Mehrani}
\ead{smehrani@fau.edu}
\author[toronto]{Samin Aref\corref{cor1}}
\ead{s.aref@utoronto.ca}
\author[gwu]{Payman Dehghanian}
\ead{payman@gwu.edu}
\cortext[cor1]{Corresponding author}
\address[toronto]{Department of Mechanical and Industrial Engineering, University of Toronto, 5 King's College Road, Toronto, Ontario M5S 3G8, Canada}
\address[fau]{Department of Information Technology and Operations Management, Florida Atlantic University, Boca Raton, Florida, USA}
\address[gwu]{Department of Electrical and Computer Engineering, The George Washington University, Washington, DC, USA}
\fi

\begin{abstract}
\singlespacing
Post-disaster power restoration requires allocating limited repair and power resources across affected communities, creating a fundamental tension between restoring service efficiently and providing equitable access to electricity. We develop an optimization framework incorporating horizontal fairness, which promotes comparable service across affected communities, and vertical fairness, which accounts for differences in community vulnerability. The framework integrates repair-crew scheduling with an approximation of power-system operations, including distributed generation and solar-plus-storage resources. It supports vulnerability-aware prioritization while protecting specified fractions of bus-level service under a horizontal-fairness reference plan and limiting the increase in restoration cost. Structural results characterize when the two fairness benchmarks can be attained simultaneously in an aggregate allocation model and provide a sufficient condition under which improving vertical fairness requires changing discrete decisions in the full restoration model. Using the IEEE 33-bus test system, we show that equitable restoration operates through two interconnected mechanisms: reconnecting isolated areas and allocating limited power within connected areas. In the tested setting, selected vulnerability-aware solutions attain both fairness benchmarks at a lower cost than their horizontal-fairness reference plan. The results also show that stronger emphasis on horizontal fairness can yield lower efficiency losses than intermediate fairness settings, with an identical-demand comparison illustrating the role of demand heterogeneity. The framework helps decision makers identify the operational sources of service disparities and evaluate how repair priorities, power allocation, and allowable cost increases affect the opportunities to reduce them.
\end{abstract}
\begin{keyword}
Fair power restoration \sep Post-disaster grid recovery \sep Equitable resource allocation \sep Horizontal fairness \sep Vertical fairness
\end{keyword}
\end{frontmatter}


\section{Introduction}
\label{s:intro}
\noindent
The global energy landscape is facing growing challenges as climate change continues to act as a threat multiplier, increasing the frequency, intensity, and duration of weather-related disruptions to power grids. Recent evidence indicates that extreme weather has become a primary driver of large-scale grid failures \citep{Panteli2015ExtremeWeather}. A 2026 report found that all of the 100 worst power-outage days in the central United States between 2014 and 2024 were weather-driven and that such extreme weather events are expected to become more frequent as the climate warms \citep{NCC2026GridResilience}. Related research has connected human-induced greenhouse gas emissions to predictable increases in extreme heat and precipitation based on thermodynamic principles and historical data \citep{Robinson2021Climate}. These events can cause extended periods of service interruption through physical damage to transmission towers, distribution lines, substations, and other components of the electric power system. Recent examples include Winter Storm Uri in 2021, which brought the Texas power grid close to a system-wide collapse \citep{Potts2024Uri}, and Super Typhoon Ragasa in 2025, which disrupted electricity service for hundreds of thousands of households in the Philippines \citep{Nando2025}. Beyond direct damage from precipitation and wind, extreme heat can reduce transmission capacity, accelerate transformer aging, and increase exposure to wildfires \citep{Gud2025TempGrid, CIEE2012California}.

Following such disruptions, system operators must decide how limited restoration resources should be allocated over time, including which damaged components should receive repair crews first and how available electricity should be distributed while parts of the network remain disconnected. These decisions are interdependent. Repairing a line changes network connectivity and may reconnect an entire group of customers to the main grid. Distributed generation and solar-plus-storage resources may also allow disconnected portions of the network to operate temporarily, although their available supply may be insufficient to satisfy all demand. Restoration planning therefore involves both the physical reconnection of affected areas and the allocation of scarce electricity as the network is progressively restored \citep{arif2018optimizing,lei2019resilient}.

Traditional restoration practices emphasize efficiency and system reliability. Utility teams generally restore the high-voltage ``backbone'' first and subsequently prioritize areas where restoration can recover substantial load within a limited amount of time. Distribution-level restoration also commonly follows a ``lifeline'' approach that gives priority to critical services essential to public safety, such as hospitals, water-treatment facilities, telecommunications, and emergency services. These principles serve important operational objectives. However, a restoration plan that performs well in terms of total restored load or aggregate outage cost does not necessarily distribute the burden of an outage fairly across affected communities. An efficiency-oriented policy may direct scarce repair resources toward locations with larger loads or greater system-wide impact while leaving smaller or geographically isolated communities without service for substantially longer periods \citep{rodriguezgarcia2024}.

The consequences of a prolonged outage also differ across communities because of differences in economic resources, demographic characteristics, access to essential services, and the ability to cope with service interruptions. The impact of an outage therefore depends on both the service shortfall and the vulnerability of those affected \citep{toplututay2024}. Assessing fairness in restoration requires attention to how service is distributed across communities and how their differing needs should influence restoration priorities.

We consider these dimensions through the notions of \emph{horizontal fairness} and \emph{vertical fairness}. Horizontal fairness promotes comparable service outcomes across affected locations. In our setting, these outcomes are measured by the average fraction of demand served at each load bus over the restoration horizon. This perspective addresses disparities that can remain hidden in aggregate restoration measures \citep{sundar2023fairly}. Vertical fairness accounts for differences in vulnerability and gives greater consideration to communities for which service interruptions may impose a greater burden. We represent vulnerability using the Social Vulnerability Index (SVI), which measures a community's vulnerability to hazardous events and its capacity to prepare for, respond to, and recover from them \citep{cutter2024SVI,mah2023social}. Our treatment of horizontal and vertical fairness builds on the resource-allocation formulation of \citet{baghersad2024tradeoffs}.

Existing studies incorporate disparities in service into power-system operation and restoration \citep{sundar2023fairly,rodriguezgarcia2024,gupta2026fair} and account for the well-being losses of vulnerable communities in grid-resilience planning \citep{toplututay2024}. \citet{sabzi2026equitable} combine equity and risk tolerance in restoration planning under uncertain line availability. Building on these developments, we examine how the interaction between horizontal and vertical fairness depends on the operational restrictions of restoration. Delayed access to a supply source limits the service a bus can receive, regardless of how power is subsequently allocated. Once loads are connected to a source, generation and storage constraints determine the scope for redistributing service. These restrictions affect whether the two fairness benchmarks can be attained together and whether an improvement requires changing discrete decisions or can be achieved through continuous operating adjustments. Understanding these distinctions is necessary for interpreting the cost of fairness and identifying the decisions that can improve it.

We develop a multi-period optimization framework for a static disaster setting in which damage, demand, and available resources are known at the beginning of the planning horizon. The system operator jointly determines repair-crew allocation, line availability, the operation of distributed generation and solar-plus-storage resources, and the fraction of demand served at each load bus. The framework uses a common set of repair and operating constraints, including an approximation of power-system operations, to compare alternative fairness preferences. The study makes three key contributions.

First, we embed horizontal and vertical fairness in a restoration framework that explicitly links vulnerability-aware prioritization to protected bus-level service and an allowable efficiency loss. A \utilitarian{} benchmark minimizes total restoration cost. The \horizontal{} model uses $\alpha$-fair utility to place progressively greater emphasis on buses with lower average service. The \horizontalVertical{} model minimizes the worst SVI-weighted cumulative service shortfall while requiring each bus to retain a specified fraction of its service under the corresponding \horizontal{} solution. When a cost limit is imposed, this reference is obtained under the same limit. This structure allows an operator to examine vulnerability-aware restoration while specifying both the service levels to be protected and the additional cost that can be tolerated. It also connects these preferences to the repair and power-allocation decisions through which service is delivered.

Second, we establish structural properties that explain how operational restrictions shape the compatibility and cost of fairness. Connectivity bounds identify service losses that cannot be recovered by reallocating power after reconnection. In an aggregate allocation benchmark with individual service ceilings, we characterize when the max-min and vertical-fairness benchmarks can be attained simultaneously. The analysis shows how a binding service ceiling can permit joint attainment, and an example demonstrates that expanding service possibilities can make the improved benchmarks incompatible. We then return to the full restoration formulation and, for fixed protected service requirements, derive a sufficient condition under which no combination of discrete line-availability, repair-crew, and battery-operating-mode decisions can be optimal at both endpoint cost limits. This result provides a criterion for distinguishing improvements that require a change in discrete decisions from those attainable through continuous operating flexibility alone.

Third, we use the IEEE 33-bus distribution test system with illustrative vulnerability scores over a 24-hour restoration horizon to examine the operational sources of fairness improvements. In the tested setting, improving service to an isolated area requires earlier restoration of a supply path, whereas vulnerability-aware service can also be improved by reallocating power among connected loads without changing the repair schedule. Selected solutions attain both the max-min and vertical-fairness benchmarks, illustrating that the two objectives do not invariably conflict. The experiments also show that a stronger emphasis on horizontal fairness can yield a lower efficiency loss than intermediate fairness settings. A comparison with identical bus-level demand helps explain the role of demand heterogeneity in this pattern, since the fairness objective evaluates service fractions while efficiency depends on the amount of demand served. These findings concern the modeled network approximation and tested configurations. They help identify when more equitable service requires changes in physical restoration and when the available power can be allocated differently.

The remainder of the paper is organized as follows. Section~\ref{s:LitRev} reviews the relevant literature on post-disaster power restoration and fairness-aware resource allocation. Section~\ref{sec:Problem} describes the post-disaster power restoration problem and the notions of horizontal and vertical fairness considered in this study. Section~\ref{sec:Formulations} presents the optimization formulations and the metrics used to evaluate the efficiency--fairness trade-offs. Section~\ref{sec:StructuralProperties} establishes structural properties of the proposed models and their fairness outcomes. Section~\ref{sec:CaseStudy} reports the numerical experiments and discusses the resulting operational and managerial insights. Section~\ref{sec:conclusions} concludes the paper and outlines directions for future research.

\section{Literature Review}
\label{s:LitRev}
\noindent
Post-disaster power restoration has been studied extensively from the perspectives of system resilience, repair scheduling, and operational efficiency. Existing studies address repair sequencing \citep{yan2020post}, integrated distribution-system restoration \citep{chen2017modernizing}, restoration costs \citep{arab2016electric}, and demand-side flexibility \citep{hafiz2019utilising}. Several models explicitly coordinate repair logistics with power-system operation. \citet{arif2018optimizing} develop a two-stage stochastic mixed-integer model that coordinates repair-crew routing, distributed generation, and network reconfiguration under uncertain repair times and demand. \citet{lei2019resilient} jointly optimize the dispatch of repair crews and mobile power sources with distribution-system restoration, accounting for transportation-network constraints and microgrid formation. \citet{li2025online} develop a simulation-based online dynamic programming method for dispatching repair crews and mobile emergency generators as information on repair times, road conditions, and failures becomes available. Related work examines restoration across interdependent infrastructure systems. \citet{maraqa2022compression} jointly determine repair sequencing and crew assignment and evaluate how additional resources can shorten restoration schedules in power, water, and gas networks. \citet{alkhaleel2022meanrisk} develop a two-stage mean-risk model that combines repair scheduling with infrastructure flows and unmet-demand costs under uncertainty, allowing multiple crews to repair a component and alternative repair modes. These studies evaluate restoration primarily through recovery time, cost, and system functionality. Network configuration also affects restoration performance. \citet{chopra2023parallel} jointly optimize network sectionalization and generator startup sequencing for parallel power-system restoration using integer programming and bounding procedures. In a related setting, \citet{geri2023corridor} develop a corridor-method matheuristic to balance customer numbers across feeders under radiality, power-flow and voltage constraints, and a limit on switching operations. Their model concerns distribution networks under normal operation, with balance measured through customer assignments across feeders.

A more recent literature explicitly incorporates fairness into power-system planning and restoration. \citet{sundar2023fairly} incorporate fairness into outage minimization through alternative objective functions and an $\epsilon$-fairness constraint, evaluating the resulting allocations using the Gini and Jain indices. \citet{fang2025energy} constrains a Grid Gini Coefficient based on load-shedding percentages to limit disparities across loads. \citet{rodriguezgarcia2024} penalize differences in the proportion of power restored across loads while allocating distributed energy resources. \citet{gupta2026fair} study reconnection after temporary disruptions and balance outage-related reliability measures and energy losses. They show that optimizing individual system objectives can produce unequal expected outage durations across industrial and residential areas. \citet{sabzi2026equitable} incorporate equity into a bi-level framework that coordinates network sectionalization and re-energization. Their model uses an approximation of the between-zone component of Theil's $T$ index to address disparities across grid sections and a gap measure for disparities within sections. A robust extension accounts for uncertain transmission-line availability, with a conservativeness parameter allowing decision makers to adjust restoration plans according to their risk tolerance.

Differences in customer needs and community vulnerability provide additional bases for evaluating restoration decisions. \citet{yan2026differentiated} incorporate differing customer reliability requirements through soft constraints in a restoration model that combines interruption and restoration costs. \citet{toplututay2024} incorporate the well-being losses of vulnerable communities into grid-resilience planning. Our treatment of fairness also draws on the broader resource-allocation literature, which examines the trade-off between aggregate efficiency and the distribution of outcomes \citep{bertsimas2012efficiency}. The $\alpha$-fairness family represents different degrees of inequality aversion. \citet{rea2021unequal} distinguish equality-oriented allocation from allocation that recognizes relevant differences among recipients, and \citet{baghersad2024tradeoffs} formalize horizontal and vertical fairness in infrastructure resource allocation by combining $\alpha$-fairness with vulnerability information. This distinction is relevant to electricity restoration because comparable service interruptions can impose different burdens across communities.

We consider horizontal and vertical fairness jointly within the operational decisions of post-disaster power restoration. The \horizontalVertical{} model incorporates vulnerability-aware priorities while protecting specified fractions of bus-level service under the corresponding \horizontal{} solution. A \utilitarian{} benchmark provides the reference for evaluating efficiency loss. These models share repair and operating constraints, with repair-crew allocation, network reconnection, and generation and solar-plus-storage dispatch determining feasible service outcomes. The structural analysis and numerical experiments examine how these restrictions affect the compatibility and cost of the two fairness objectives. In particular, we distinguish improvements that require changes in discrete restoration decisions from those attainable by reallocating power among connected loads.

\section{Post-disaster Power Restoration Problem}
\label{sec:Problem}
\noindent
We consider post-disaster restoration of a distribution system over a discrete planning horizon $\mathcal{T}$. We focus on a static restoration setting, in which the damage caused by the disaster is observed at the beginning of the planning horizon. Accordingly, damaged power lines, load demands, available repair crew members, and local generation and storage resources are treated as known inputs. Although this information is fixed at the outset, restoration decisions are made over multiple periods as damaged lines are repaired and available power is dispatched throughout the system.

In each period, the system operator determines the allocation of repair crew members to damaged lines, the dispatch of available generation and solar-plus-storage resources, and the fraction of demand served at each load bus. The pre-disaster system is assumed to have sufficient generation capacity to meet demand under normal operating conditions. Post-disaster service interruptions therefore result from damaged network components and the limited resources available during restoration. We measure system efficiency by the fraction of total demand served over the restoration horizon. Restoration decisions, however, also determine how this service is distributed across affected buses. We therefore consider horizontal fairness based on the balance of service outcomes across load buses and vertical fairness based on the vulnerability of the communities they serve. For horizontal fairness, the service outcome of each load bus is measured by its average fraction of demand served over the restoration horizon. Vertical fairness accounts for community vulnerability through SVI, giving greater consideration to unmet demand at buses serving more vulnerable communities. The resulting problem captures the trade-off between restoring as much load as possible, balancing service across affected buses, and accounting for differences in vulnerability.

The power distribution system is represented by a set of buses $\mathcal{B}$ connected by a set of distribution lines $\mathcal{E}$. The complete notation used throughout the model is provided in Section~\ref{sec:NOMENCLATURE}. Bus 1 denotes the substation bus through which the distribution system receives power from the transmission system, and $\mathcal{B}^{L}=\mathcal{B}\setminus\{1\}$ denotes the set of load buses. A subset $\mathcal{B}^{c}\subseteq\mathcal{B}^{L}$ serves critical loads, such as hospitals and emergency facilities, whose loss of service may have significant health, safety, or emergency-response consequences. The system also includes a set $\mathcal{G}$ of controllable distributed generators and a set $\mathcal{S}$ of solar-plus-storage systems. Each solar-plus-storage system consists of photovoltaic (PV) panels, battery storage, and an inverter. For each bus $i\in\mathcal{B}$, $\mathcal{G}_i$ and $\mathcal{S}_i$ denote the controllable generators and solar-plus-storage systems located at that bus, respectively.

At the beginning of the planning horizon ($t=0$), a known subset of distribution lines $\mathcal{E}^{f}\subseteq\mathcal{E}$ is out of service due to disaster-related damage. Starting in period $t=1$, up to $\overline{C}$ repair crew members can be assigned in each period to restore the failed power lines. Multiple crew members may work on the same line, and each failed line $\ell\in\mathcal{E}^{f}$ requires $C_{\ell}$ crew-hours to restore. Crew members assigned to an unfinished line remain on that line until its repair is completed. Once fixed, the line becomes available for the remainder of the planning horizon, and its assigned crew members may be reassigned in the following period.

\subsection{Nomenclature}\label{sec:NOMENCLATURE}
\noindent
The notation used to formulate the restoration decisions and the associated network operations is summarized as follows.

\subsubsection{Sets}
\begingroup
\singlespacing\small
\begin{longtable*}{@{}>{\raggedright\arraybackslash}p{0.31\linewidth}>{\raggedright\arraybackslash}p{\dimexpr0.69\linewidth-2\tabcolsep\relax}@{}}
$\mathcal{T}$ & Set of restoration periods, indexed by $t \in \mathcal{T}=\{1,\ldots,\tau\}$. \\
$\mathcal{T}_{\ge 2}$ &  Set of restoration periods after the first period, $\mathcal{T}_{\ge 2}=\{2,\ldots,\tau\}$. \\
$\mathcal{B}$ & Set of buses, indexed by $i,j\in\mathcal{B}$. Bus 1 is the substation bus connecting to the transmission system. \\
$\mathcal{B}^L$ &Set of load buses, $\mathcal{B}^{L}=\mathcal{B}\setminus\{1\}$. \\
$\mathcal{B}^c$ &Set of buses serving critical loads, with $\mathcal{B}^c \subseteq \mathcal{B}^L$. \\
$\mathcal{B}_i$ & Set of buses adjacent to bus $i$. \\
$\mathcal{E}$ & Set of distribution lines, indexed by $\ell=(i,j)\in\mathcal{E}$. \\
$\mathcal{E}^f$ & Set of lines failed immediately after the disaster, $\mathcal{E}^f \subseteq \mathcal{E}$. \\
$\mathcal{G}$ & Set of controllable distributed generators, indexed by
$g\in \mathcal{G}$. \\
$\mathcal{S}$ & Set of solar-plus-storage systems, indexed by $s\in \mathcal{S}$. \\
$\mathcal{S}_i$ & Subset of solar-plus-storage systems located at bus $i$, with $\mathcal{S}_i\subseteq\mathcal{S}$. \\
$\mathcal{G}_i$ & Subset of controllable generators located at bus $i$, with $\mathcal{G}_i\subseteq \mathcal{G}$. \\
\end{longtable*}
\endgroup

\subsubsection{Parameters}

\begingroup
\singlespacing\small
\begin{longtable*}{@{}>{\raggedright\arraybackslash}p{0.31\linewidth}>{\raggedright\arraybackslash}p{\dimexpr0.69\linewidth-2\tabcolsep\relax}@{}}
$\Delta t$ & Length of each restoration period (in hours). \\

$p_{i,t}^D$ & Active power demand at bus $i$ in period $t$. \\
$q_{i,t}^D$ & Reactive power demand at bus $i$ in period $t$. \\
$g_{\ell}$ & Conductance of line $\ell$. \\
$b_{\ell}$ & Susceptance of line $\ell$. \\
$b^{sh}_{\ell}$ & Shunt susceptance of line $\ell$. \\
$\overline{S}_{\ell}$ & Apparent power capacity of line $\ell$. \\
$\overline{v}_i$ & Voltage upper limit at bus $i$. \\
$a_{\ell,0}$ & Binary parameter equal to 1 if line $\ell$ is available immediately after the disaster. \\
$\overline{p}^{S}_{s,t}$ & Maximum available active power from solar system $s$ in period $t$. \\
$\overline{S}^{I}_{s}$ & Apparent power capacity of smart inverter $s$. \\
$\underline{PF}_{s}$ & Minimum power factor of smart inverter $s$. \\
$\overline{p}^{B}_{s}$ & Maximum charging or discharging rate of battery system $s$. \\
$\eta^{in}_{s}$ & Charging efficiency of battery system $s$. \\
$\eta^{out}_{s}$ & Discharging efficiency of battery system $s$. \\
$\underline{E}_{s}, \overline{E}_{s}$ &  Minimum and maximum stored energy in battery system $s$. \\
$E^0_{s}$ & Initial energy stored in battery system $s$. \\
$\overline{p}^{G}_{g}$ & Maximum active power output of controllable generator $g$. \\
$\overline{S}^{G}_{g}$ & Apparent power capacity of controllable generator $g$. \\
$\overline{F}_{g}$ & Available energy of controllable generator $g$ over the restoration horizon. \\
$\overline{C}$ & Total number of available repair crew members in each period. \\
$C_{\ell}$ & Required repair crew-hours to restore failed line $\ell\in\mathcal{E}^{f}$. \\
$c^{gen}$ & Unit cost of supplied energy. \\
$VOLL_i$ & Unit value of lost load at bus $i$. \\
$SVI_i$ &SVI associated with bus $i$. \\
$\alpha$ & Horizontal-fairness parameter in the horizontal fairness formulation. \\
$\beta$ & Vertical-fairness parameter in the horizontal-vertical-fairness formulation. \\
$\lambda$ & Budget factor limiting the allowable efficiency loss under fairness-aware restoration. \\
$Z^{U}$ &Optimal objective value of the baseline utilitarian restoration model. \\
\end{longtable*}
\endgroup

\subsubsection{Decision Variables}

\begingroup
\singlespacing\small
\begin{longtable*}{@{}>{\raggedright\arraybackslash}p{0.31\linewidth}>{\raggedright\arraybackslash}p{\dimexpr0.69\linewidth-2\tabcolsep\relax}@{}}
$\gamma_{i,t}\in[0,1]$ & Fraction of demand served at load bus $i\in\mathcal{B}^L$ in period $t$. \\
$y_i \in [0,1]$ & Average fraction of power demand served at load bus $i$ over the restoration horizon. \\
$a_{\ell,t}\in\{0,1\}$ & Binary variable equal to 1 if line $\ell\in\mathcal{E}$ is available in period $t$. \\
$p^G_{1,t}\in\mathbb{R}_{+},\ q^G_{1,t}\in\mathbb{R}$ & Active and reactive power received from the transmission system in period $t$. \\
$p^S_{s,t} \in \mathbb{R}_{+}$ & Active solar power from system $s$ directly supplied from panels to the grid in period $t$. \\
$E_{s,t} \in [\underline{E}_{s}, \overline{E}_{s}]$ & Energy stored in battery system $s$ at the end of period $t$. For $t=0$, $E_{s,0}=E^0_{s}$. \\
$c^B_{s,t} \in  \mathbb{R}_{+}$ & Charging rate of battery system $s$ in period $t$. \\
$d^B_{s,t} \in \mathbb{R}_{+}$ & Discharging rate of battery system $s$ in period $t$. \\
$z^B_{s,t}\in\{0,1\}$ & Charging-mode indicator for battery system $s$ in period $t$. \\
$q^I_{s,t} \in \mathbb{R}$ & Reactive power supplied or absorbed by inverter $s$ in period $t$. \\
$p^G_{g,t}, q^G_{g,t} \in \mathbb{R}$ & Active and reactive power output of controllable generator $g$ in period $t$. \\
$p_{t}^\ell, q_{t}^\ell \in \mathbb{R}$ & Active and reactive power flow on line $\ell$ in period $t$. \\
$u_{i,\ell,t} \in \mathbb{R}_{+}$ & Auxiliary voltage variable associated with bus $i$ at one end of line $\ell$ in period $t$. \\
$R_{\ell,t}, T_{\ell,t} \in \mathbb{R}$ & Auxiliary variables used in the AC power-flow approximation for line $\ell$ in period $t$. \\
$v_{i,t} \in \mathbb{R}_{+}$ & Squared voltage magnitude at bus $i$ in period $t$, divided by $\sqrt{2}$. \\
$r_{\ell,t}\in\mathbb{Z}_+$ & Number of repair crew members assigned to failed line $\ell\in\mathcal{E}^f$ in period $t$. \\
\end{longtable*}
\endgroup

\setcounter{table}{0}
\renewcommand{\theHtable}{main.\arabic{table}}

\section{Problem Formulations}
\label{sec:Formulations}
\noindent
We formulate the post-disaster power restoration problem described in Section~\ref{sec:Problem} by first defining the feasible set of restoration plans, denoted by $\Omega$, through a common set of operational and network constraints. These constraints capture demand service, power balance, approximate power-flow relationships, voltage and line-status limits, distributed generation, solar-plus-storage operation, and repair-crew scheduling. We then introduce the \utilitarian{} benchmark, which selects an efficient restoration plan from $\Omega$ by minimizing total restoration cost, including generation cost and the value of lost load. Building on the same feasible set, we present two fairness-aware formulations: the \horizontal{} model, based on the average fraction of demand served at each load bus, and the \horizontalVertical{} model, which further accounts for community vulnerability through SVI-weighted unmet demand. Finally, we define the trade-off metrics used to compare the models in terms of efficiency, horizontal fairness, and vertical fairness.

\subsection{Feasible Restoration Set}
\label{sec:FeasibleRestorationSet}
\noindent
The feasible restoration set is defined by Constraints~\eqref{eq:power_balance}--\eqref{eq:repair}, which capture power balance and power-flow requirements, network availability, the operation of distributed generation and solar-plus-storage systems, and the allocation of repair crews to damaged lines. The power-system constraints are adapted from \citet{rodriguezgarcia2024}, with modifications to represent the operation of solar-plus-storage systems. Each such system consists of PV panels, dedicated battery storage, and an inverter that controls its interaction with the grid. All constraints listed below are defined for all $t \in \mathcal{T}$ unless noted otherwise.

\subsubsection{Power Balance and Basic Limits}
\noindent
Constraint~\eqref{eq:1_power_balance} defines the active and reactive power received from the transmission system at the substation based on the corresponding outgoing power flows. Constraints~\eqref{eq:bus_p_balance} and \eqref{eq:bus_q_balance} enforce active and reactive power balance, respectively, at each load bus. Specifically, available power from controllable generation and solar-plus-storage systems is balanced against served demand and power flows through the power distribution lines. These balances neglect branch losses, with each line flow defined in a fixed orientation. Active power imports are nonnegative, with no separate import cap beyond the network limits. Reactive power exchange may have either sign.

\begin{subequations}
\label{eq:power_balance}
\begin{align}
    p_{1,t}^G = \sum_{j \in \mathcal{B}_1} p^{(1,j)}_t, &\quad q_{1,t}^G = \sum_{j \in \mathcal{B}_1} q^{(1,j)}_t,  \label{eq:1_power_balance}\\
     \sum_{s \in \mathcal{S}_i} (p_{s,t}^S + d_{s,t}^B\eta^{out}_{s})+ \sum_{g \in \mathcal{G}_i} p_{g,t}^G  & = \gamma_{i,t} \, p_{i,t}^D + \sum_{\ell=(i,j)\in\mathcal{E}}p^\ell_t-\sum_{\ell=(j,i)\in\mathcal{E}}p^\ell_t, \quad \forall i \in \mathcal{B}^L \label{eq:bus_p_balance} \\
    \sum_{s \in \mathcal{S}_i} q_{s,t}^I + \sum_{g \in \mathcal{G}_i} q_{g,t}^G& = \gamma_{i,t}\, q_{i,t}^D + \sum_{\ell=(i,j)\in\mathcal{E}}q^\ell_t-\sum_{\ell=(j,i)\in\mathcal{E}}q^\ell_t, \quad \forall i \in \mathcal{B}^L \label{eq:bus_q_balance}
\end{align}
\end{subequations}

\subsubsection{Power Flow Constraints}
\noindent
Constraints~\eqref{eq:power_flow} approximate the power-flow relationships in the distribution network. The first two constraints relate active and reactive power flows to bus-voltage and line parameters, while the third limits the apparent power flow on each line by its capacity. The next four constraints link the auxiliary voltage variables to line availability, ensuring that voltage relationships are enforced only when a line is in service. The final constraint provides the convex relaxation relating the auxiliary voltage variables used in the network approximation. The substation voltage is fixed by $v_{1,t}=\frac{(\overline{v}_1)^2}{\sqrt{2}}$. The optimization imposes nonnegativity and upper voltage bounds, without a positive lower voltage bound. The reported service and fairness outcomes are evaluated under this network approximation and do not establish feasibility under the full AC power-flow equations.

\begin{subequations}
\label{eq:power_flow}
\begin{align}
p_t^{\ell} &= \sqrt{2} g_{\ell}\, u_{i,\ell,t} - g_{\ell} \, R_{\ell,t} - b_{\ell}\, T_{\ell,t}, \quad &&\forall \ell=(i,j)\\
q_t^{\ell} &= -\sqrt{2} \left( b_{\ell} + \frac{b_{\ell}^{sh}}{2} \right) u_{i,\ell,t} + b_{\ell} R_{\ell,t} - g_{\ell}\, T_{\ell,t}, \quad &&\forall \ell=(i,j)\\
(p_t^{\ell})^2 &+ (q_t^{\ell})^2 \leq (\overline{S_{\ell}})^2, \quad &&\forall \ell=(i,j)\\
0 &\leq u_{i,\ell,t} \leq \frac{\overline{v_i}^2}{\sqrt{2}} a_{\ell,t}, \quad &&\forall \ell=(i,j)\\
0 &\leq u_{j,\ell,t} \leq \frac{\overline{v_j}^2}{\sqrt{2}} a_{\ell,t}, \quad &&\forall \ell=(i,j)\\
0 &\leq v_{i,t} - u_{i,\ell,t} \leq \frac{\overline{v_i}^2}{\sqrt{2}} (1 - a_{\ell,t}), \quad &&\forall \ell=(i,j)\\
0 &\leq v_{j,t} - u_{j,\ell,t} \leq \frac{\overline{v_j}^2}{\sqrt{2}} (1 - a_{\ell,t}), \quad &&\forall \ell=(i,j) \\
2 &u_{i,\ell,t}\,u_{j,\ell,t} \geq (R_{\ell,t})^2 + (T_{\ell,t})^2, \quad &&\forall \ell=(i,j)
\end{align}
\end{subequations}

\subsubsection{Solar-Plus-Storage System Constraints}
\noindent
PV generation can be supplied directly to the grid or used to charge the battery, while stored energy can be discharged to support the grid when needed. Constraint~\eqref{eq:b_charge} tracks the battery state of charge based on its charging and discharging decisions, and Constraint~\eqref{eq:charge_limit} limits both rates by the battery power capacity and prevents simultaneous charging and discharging. Constraint~\eqref{eq:inverter_limit} ensures that the power supplied directly from the PV panels and the power used for battery charging do not exceed the available PV generation, accounting for charging efficiency. Constraint~\eqref{eq:inverter_capacity} limits the combined active and reactive power supplied through the inverter by its apparent power capacity. Finally, Constraints~\eqref{eq:inverter_PF1} and \eqref{eq:inverter_PF2} bound reactive power supply and absorption according to the inverter's minimum power factor.

\begin{subequations}
\label{eq:solar}
\begin{align}
& E_{s, t} = E_{s,t-1} + (-d_{s,t}^B + c_{s,t}^B)\Delta t, \quad &&\forall s \label{eq:b_charge}\\
& 0\leq c_{s,t}^B\leq\overline{p_s^B}z^B_{s,t},\quad 0\leq d_{s,t}^B\leq\overline{p_s^B}(1-z^B_{s,t}), \quad &&\forall s, \label{eq:charge_limit} \\
& p_{s,t}^S + c_{s,t}^B/\eta_{s}^{in} \leq \overline{p_{s,t}^S}, \quad &&\forall s, \label{eq:inverter_limit}\\
& \left( q_{s,t}^I \right)^2 \le \left( \overline{S_s^I} \right)^2 - \left( p_{s,t}^S + d_{s,t}^B\eta_{s}^{out} \right)^2, \quad &&\forall s, \label{eq:inverter_capacity}\\
& q_{s,t}^I \le \left( p_{s,t}^S + d_{s,t}^B\eta_{s}^{out}\right) \sqrt{1/\underline{PF_s}^2 - 1}, \quad &&\forall s, \label{eq:inverter_PF1} \\
& q_{s,t}^I \ge -\left( p_{s,t}^S + d_{s,t}^B\eta_{s}^{out}\right) \sqrt{1/\underline{PF_s}^2 - 1}, \quad &&\forall s. \label{eq:inverter_PF2}
\end{align}
\end{subequations}

\subsubsection{Distributed Controllable Generator Constraints}
\noindent
Constraint~\eqref{eq:cg_p} limits the active power output of each controllable generator in each period, while Constraint~\eqref{eq:cg_E} limits its total energy production over the restoration horizon. Constraint~\eqref{eq:cg_total} further restricts the combined active and reactive power output by the generator's apparent power capacity.

\begin{subequations}
\label{eq:generator}
\begin{align}
    &0 \le p_{g,t}^G \le \overline{p_{g}^G}, && \quad \forall g, \label{eq:cg_p}\\
   & \sum_{t \in \mathcal{T}} p_{g,t}^G \Delta t \le \overline{F_{g}}, && \quad \forall g \label{eq:cg_E}\\
   &\left( q_{g,t}^G \right)^2 + \left( p_{g,t}^G \right)^2 \le \left( \overline{S_{g}^G} \right)^2 , && \quad \forall g. \label{eq:cg_total}
\end{align}
\end{subequations}

\subsubsection{Power Line Outage and Repair Constraints}
\noindent
Constraints~\eqref{eq:max_crew}--\eqref{eq:crew_front_load} govern the allocation of repair crews and the restoration of failed lines. Constraint~\eqref{eq:max_crew} limits the total number of crew members assigned in each period by the available crew capacity $\overline{C}$. Constraint~\eqref{eq:req_work} allows a failed line to become available only after its required repair crew-hours have been completed, while Constraint~\eqref{eq:line_status_monotonic} ensures that a restored line remains available thereafter. Constraint~\eqref{eq:crew_monotonic} requires the number of crew members assigned to an unfinished line to remain at least at its previous-period level, and Constraint~\eqref{eq:no_crew_if_fixed} prevents further crew assignment once the line has been restored. Together, these two constraints prevent crew members from being withdrawn from an unfinished repair and reassigned elsewhere. Finally, Constraint~\eqref{eq:crew_front_load} requires the total number of deployed crew members to be nonincreasing over the restoration horizon, so that deployed crew capacity cannot increase after the first period.

\begin{subequations}
\label{eq:repair}
\begin{align}
& \sum_{\ell \in \mathcal{E}^f} r_{\ell,t} \le \overline{C}, \quad &&\label{eq:max_crew} \\
& C_{\ell} \cdot a_{\ell,t} \le \Delta t\sum_{\tau=1}^{t-1} r_{\ell,\tau}, \quad && \forall \ell \in \mathcal{E}^f \label{eq:req_work} \\
& a_{\ell,t} \ge a_{\ell,t-1}, \quad && \forall \ell \in \mathcal{E}^f \label{eq:line_status_monotonic} \\
& r_{\ell,t} \ge r_{\ell,t-1} - \overline{C} \cdot a_{\ell,t}, \quad && \forall \ell \in \mathcal{E}^f \label{eq:crew_monotonic} \\
& r_{\ell,t} \le \overline{C} \cdot (1 - a_{\ell,t}), \quad &&  \forall \ell \in \mathcal{E}^f \label{eq:no_crew_if_fixed} \\
& \sum_{\ell \in \mathcal{E}^f} r_{\ell,t} \leq \sum_{\ell \in \mathcal{E}^f} r_{\ell,t-1}, \quad &&  \forall t \in \mathcal{T}_{\ge 2} \label{eq:crew_front_load}
\end{align}
\end{subequations}

We set $r_{\ell,0}=0$ and $a_{\ell,0}=a_{\ell,1}=0$ for failed lines, while healthy lines remain available throughout the horizon. A line can be used in period $t$ only if its required crew-hours have been completed by the end of period $t-1$. Its crews may then be reassigned in period $t$.

\subsection{Utilitarian Benchmark}
\noindent
The \utilitarian{} model serves as the efficiency benchmark by selecting, among the feasible restoration plans, the plan that minimizes the total cost over the restoration horizon. As shown in \eqref{eq:ut_obj}, this cost consists of the cost of supplied electricity and the cost associated with unserved demand, measured using VOLL. Because VOLL represents the economic consequence of service interruptions and is typically substantially higher than the cost of supplying electricity, the objective favors restoring and serving demand as early as possible. The \utilitarian{} model is formulated as

\begin{equation}
\label{eq:ut_obj}
    \min \sum_{t \in \mathcal{T}} \sum_{i \in \mathcal{B}^L} c^{gen}\gamma_{i,t}p^D_{i,t}\Delta t + VOLL_i(1-\gamma_{i,t})p^D_{i,t}\Delta t
\end{equation}

\subsection{Horizontal-Fairness Formulation}
\label{sec:HE}
\noindent
We model horizontal fairness using the $\alpha$-fairness framework \citep{ogryczak2014fair,baghersad2024tradeoffs}. For an outcome $y_i>0$, the corresponding utility function is
\begin{align}
\label{eq:W}
    W_{\alpha}(y_i) =
\begin{cases}
\frac{y_i^{1-\alpha}}{1-\alpha}, & \text{if } \alpha \neq 1, \\
\log(y_i), & \text{if } \alpha = 1.
\end{cases}
\end{align}
The parameter $\alpha>0$ controls the emphasis placed on balancing outcomes across recipients, with larger values placing greater weight on those receiving lower service levels.

In our restoration setting, the outcome for each load bus is measured by its average fraction of demand served over the restoration horizon:
\begin{align}
    y_i = \frac{1}{|\mathcal{T}|}\sum_{t \in \mathcal{T}}\gamma_{i,t},
    \quad \forall i\in \mathcal{B}^L.
    \label{eq:yi}
\end{align}
Using a fraction rather than the total amount of energy supplied allows service outcomes to be compared across buses with different demand levels. Let $\Gamma_0$ denote the set of restoration decisions satisfying the feasible restoration set defined in Section~\ref{sec:FeasibleRestorationSet}. The \horizontal{} model is then formulated as
\begin{align}
    W_{\alpha} = \max_{\gamma \in \Gamma_0}
\begin{cases}
\displaystyle \sum_{i \in \mathcal{B}^L}
\frac{\left[\frac{1}{|\mathcal{T}|}\sum_{t \in \mathcal{T}}\gamma_{i,t}\right]^{1-\alpha}}
{1-\alpha},
& \text{if } \alpha>0,\ \alpha\neq 1, \\[2mm]
\displaystyle \sum_{i \in \mathcal{B}^L}
\log\left(\frac{1}{|\mathcal{T}|}\sum_{t \in \mathcal{T}}\gamma_{i,t}\right),
& \text{if } \alpha=1.
\end{cases}
\label{eq:af_2}
\end{align}
When $\alpha=1$, \eqref{eq:af_2} maximizes the sum of logarithms of the average service fractions, corresponding to the objective associated with Proportional Fairness (PF) \citep{Kelly1997PF}. Its marginal utility, $1/y_i$, assigns greater weight to improvements at buses with lower average service levels. On a convex feasible set of strictly positive service vectors, maximizing this objective is equivalent to requiring that every feasible alternative yield a non-positive aggregate proportional change relative to the optimal allocation. In our restoration model, discrete line-availability and repair-crew decisions generally make the feasible set nonconvex. We therefore use PF to denote the logarithmic objective at $\alpha=1$, without asserting that its optimizer satisfies this proportional-change property across all feasible restoration plans.
As $\alpha$ increases, the objective places progressively greater emphasis on buses with lower service levels. In the limiting case as $\alpha\rightarrow\infty$, the \horizontal{} model approaches Max-Min Fairness (MMF), which maximizes the service level of the worst-served bus:
\begin{align}
    \max_{\gamma \in \Gamma_0}
    \left\{
    \min_{i\in \mathcal{B}^L}
    \frac{1}{|\mathcal{T}|}\sum_{t \in \mathcal{T}}\gamma_{i,t}
    \right\}.
    \label{eq:mmf_obj}
\end{align}
The unconstrained \horizontal{} model characterizes the restoration plan associated with a given fairness preference without restricting the resulting increase in restoration cost relative to the \utilitarian{} benchmark. To examine horizontal fairness under a limited efficiency loss, we additionally impose
\begin{align}
    \sum_{t \in \mathcal{T}} \sum_{i \in \mathcal{B}^L}
    \left[
    c^{gen}\gamma_{i,t}p^D_{i,t}\Delta t
    + VOLL_i(1-\gamma_{i,t})p^D_{i,t}\Delta t
    \right]
    \leq \lambda Z^U,
    \label{eq:budget_con}
\end{align}
where $Z^U$ is the optimal objective value of the \utilitarian{} benchmark and $\lambda\geq1$ specifies the allowable increase in total restoration cost. Thus, $\lambda=1$ restricts the \horizontal{} model to restoration plans with the same minimum cost as the \utilitarian{} benchmark, while larger values allow additional cost in pursuit of horizontal fairness.
\subsection{Horizontal-Vertical-Fairness Formulation}
\noindent
The \horizontalVertical{} model extends the \horizontal{} model by accounting for differences in vulnerability across load buses. Following \citet{baghersad2024tradeoffs}, we use the \horizontal{} solution as a service-level benchmark and allow departures from this benchmark to improve vertical fairness. Let $y_i^{\alpha}$ denote the average fraction of demand served at bus $i$ under the \horizontal{} solution for a given $\alpha$. The parameter $\beta\in[0,1]$ controls the extent to which these service levels must be preserved. The \horizontalVertical{} model is formulated as
\begin{subequations}
\label{eq:hve}
\begin{align}
& \min_{\mathbf{\gamma} \in \Gamma_0} \left\{ \max_{i \in \mathcal{B}^L} \quad
 SVI_i \sum_{t \in \mathcal{T}} (1 - \gamma_{i,t}) \right\},
\label{eq:hve_obj} \\
& \quad \quad \text{s.t.} \quad  y_i \geq \beta \, y_i^{\alpha}, \quad \forall i \in \mathcal{B}^L.
\label{eq:hve_cst}
\end{align}
\end{subequations}
The objective in \eqref{eq:hve_obj} minimizes the largest SVI-weighted accumulated fraction of unmet demand across load buses, thereby placing greater emphasis on avoiding prolonged service shortfalls at more vulnerable locations. Constraint~\eqref{eq:hve_cst} requires each bus to retain at least a fraction $\beta$ of the average service it receives under the corresponding \horizontal{} solution. When $\beta=1$, these service levels must be fully preserved, while decreasing $\beta$ provides greater flexibility to reallocate service based on vulnerability. At $\beta=0$, the service-level requirement becomes nonbinding, and restoration decisions are driven by the vertical-fairness objective subject to the common feasible restoration set.
When examining the \horizontalVertical{} model under a limited efficiency loss, we additionally impose the budget constraint in \eqref{eq:budget_con}.
For a finite cost factor $\lambda$, let $\mathbf{x}^{HF,\alpha,\lambda}$ denote the horizontal-fairness solution under the same cost limit, with average service vector $\mathbf{y}^{\alpha,\lambda}$. The service-level requirement in \eqref{eq:hve_cst} is then replaced by
\begin{equation}
 y_i\geq\beta y_i^{\alpha,\lambda},\qquad \forall i\in\mathcal{B}^L.
 \label{eq:budget_service_floor}
\end{equation}
For $0\leq\beta\leq1$, this horizontal reference is feasible for the corresponding horizontal-vertical-fairness model. Hence,
\begin{equation}
 V^{HVF,*}_{\alpha,\beta,\lambda}\leq V\!\left(\mathbf{x}^{HF,\alpha,\lambda}\right).
 \label{eq:hf_reference_dominance}
\end{equation}
The same ordering holds for VFD under a common benchmark and positive normalization denominator.
\subsection{Trade-off Metrics}
\noindent
The three formulations reflect different priorities in post-disaster restoration, making it useful to quantify the trade-offs among efficiency, horizontal fairness, and vertical fairness. Following the price-of-fairness framework \citep{bertsimas2011price,bertsimas2012efficiency}, we use the Price of Fairness (POF) to measure the efficiency loss associated with fairness-aware restoration. We then measure departures from the two fairness benchmarks using Deviation from Max-Min Fairness (DMMF) and Vertical-Fairness Deviation (VFD). The former evaluates the worst-bus service outcome relative to MMF, while the latter evaluates the vulnerability-aware outcome relative to its best attainable value. Finally, we complement these benchmark-relative measures with Jain-based indices that characterize the overall distribution of horizontal- and vertical-fairness outcomes across load buses.
\subsubsection{Price of Fairness}
\noindent
POF measures the loss in system efficiency incurred by a fairness-aware allocation relative to the \utilitarian{} benchmark \citep{bertsimas2011price,bertsimas2012efficiency}. In our setting, we measure efficiency by the total demand served over the restoration horizon. Accordingly, for a fairness-aware model $F$, we define
\begin{equation}
    \text{POF}(F) =
    \frac{
    \displaystyle
    \sum_{i\in\mathcal{B}^L}\sum_{t\in\mathcal{T}}
    \left(\gamma_{i,t}p^D_{i,t}\right)^U
    -
    \displaystyle
    \sum_{i\in\mathcal{B}^L}\sum_{t\in\mathcal{T}}
    \left(\gamma_{i,t}p^D_{i,t}\right)^F
    }{
    \displaystyle
    \sum_{i\in\mathcal{B}^L}\sum_{t\in\mathcal{T}}
    \left(\gamma_{i,t}p^D_{i,t}\right)^U
    },
    \label{eq:poff2}
\end{equation}
where the superscripts $U$ and $F$ denote the solutions obtained from the \utilitarian{} and fairness-aware models, respectively. Larger values of POF indicate a greater reduction in total demand served relative to the \utilitarian{} benchmark.
\subsubsection{Deviation from Max-Min Fairness}
\noindent
To assess horizontal fairness, we compare the worst-bus service outcome of a restoration plan with that achieved under MMF. Because MMF maximizes the minimum average fraction of demand served across load buses, it provides a natural benchmark for this comparison. We define DMMF as
\begin{equation}
    \text{DMMF} =
    \frac{
    \displaystyle
    \min_{i\in\mathcal{B}^L} y_i^{MMF}
    -
    \displaystyle
    \min_{i\in\mathcal{B}^L} y_i
    }{
    \displaystyle
    \min_{i\in\mathcal{B}^L} y_i^{MMF}
    },
    \label{eq:dmmf}
\end{equation}
where $y_i^{MMF}$ denotes the average fraction of demand served at bus $i$ under the MMF solution. DMMF therefore represents the relative shortfall in the service level of the worst-served bus from the highest attainable minimum service level. A value of zero indicates that the restoration plan attains the same minimum service level as MMF, although service outcomes for the remaining buses may differ.
\subsubsection{Vertical-Fairness Deviation}
\noindent
Vertical fairness is evaluated relative to the best attainable outcome of the vulnerability-aware objective. Let
\begin{equation}
    V =
    \max_{i\in\mathcal{B}^L}
    \left\{
    SVI_i \sum_{t\in\mathcal{T}}(1-\gamma_{i,t})
    \right\}
    \label{eq:vf_outcome}
\end{equation}
denote the largest SVI-weighted accumulated fraction of unmet demand under a given restoration plan. We use $V^{VF}$ to denote the minimum attainable value of this measure when the service-level requirement in \eqref{eq:hve_cst} is nonbinding and no efficiency limit is imposed. We define VFD as
\begin{equation}
    \text{VFD} =
    \frac{
    V-V^{VF}
    }{
    \displaystyle
    \max_{i\in\mathcal{B}^L}\{SVI_i\}|\mathcal{T}|-V^{VF}
    }.
    \label{eq:vfd}
\end{equation}
VFD measures the deterioration in the vulnerability-aware outcome relative to its best attainable value, normalized by the range between that benchmark and the maximum possible SVI-weighted unmet-service outcome. A value of zero indicates that the restoration plan attains the vertical-fairness benchmark, while larger values indicate a greater departure from it.
\subsubsection{Fairness Indices}
\noindent
DMMF and VFD evaluate restoration plans relative to benchmarks based on the worst-served bus. To complement these measures, we use distributional indices that summarize service outcomes across all load buses. For horizontal fairness, we adopt the Jain Index \citep{jain1984fairness}:
\begin{equation}
    \text{JI} =
    \frac{
    \left(\sum_{i\in\mathcal{B}^L} y_i\right)^2
    }{
    |\mathcal{B}^L|
    \sum_{i\in\mathcal{B}^L} y_i^2
    }.
    \label{eq:jain}
\end{equation}
JI ranges from $1/|\mathcal{B}^L|$ to $1$, with larger values indicating a more balanced distribution of average service across load buses. In particular, $\text{JI}=1$ when all load buses receive the same average fraction of demand served over the restoration horizon.
To provide a corresponding system-wide measure of vertical fairness, we define a vulnerability-adjusted Jain Index based on the SVI-weighted unmet-service outcomes used in the \horizontalVertical{} objective:
\begin{equation}
    \text{VJI} =
    \frac{
    \left(\sum_{i\in\mathcal{B}^L} SVI_i(1-y_i)\right)^2
    }{
    |\mathcal{B}^L|
    \sum_{i\in\mathcal{B}^L}
    \left(SVI_i(1-y_i)\right)^2
    }.
    \label{eq:vjain}
\end{equation}
VJI likewise ranges from $1/|\mathcal{B}^L|$ to $1$ whenever some unmet demand remains. Larger values indicate a more balanced distribution of vulnerability-weighted unmet service across load buses, with $\text{VJI}=1$ when $SVI_i(1-y_i)$ is identical across all buses. If all buses are fully served, VJI is defined as $1$.
\section{Structural Properties}
\label{sec:StructuralProperties}
\noindent
The formulations developed in the previous section evaluate fairness through service accumulated over a multi-period restoration process. The timing of this service is constrained by network connectivity, repair decisions, and the availability of local power resources. We next establish properties that clarify how these operational restrictions shape attainable fairness outcomes. We first identify when different restoration trajectories are equivalent from a fairness perspective and how delayed access to power limits service. We then characterize when horizontal and vertical fairness can be attained simultaneously in an aggregate allocation benchmark. Finally, we return to the full formulation and distinguish improvements that can be achieved through continuous operating decisions from those that require a change in discrete restoration decisions. Proofs are provided in \ref{app:proofs}.
Both the \horizontal{} and \horizontalVertical{} formulations evaluate a bus using its service accumulated over the restoration horizon. Consequently, different restoration trajectories can produce the same fairness outcomes if they deliver the same cumulative service to every bus.
\begin{lemma}[Temporal invariance of fairness outcomes]
\label{lem:temporal_invariance}
Consider two feasible restoration plans $\boldsymbol{\gamma}$ and $\tilde{\boldsymbol{\gamma}}$ such that
\[
\sum_{t\in\mathcal{T}}\gamma_{i,t}
=
\sum_{t\in\mathcal{T}}\tilde{\gamma}_{i,t},
\qquad
\forall i\in\mathcal{B}^L.
\]
The two plans have identical \horizontal{} objective values for any $\alpha$ and identical \horizontalVertical{} objective values. They also have identical DMMF, JI, VFD, and VJI, whenever these quantities are defined and evaluated using the same benchmarks.
\end{lemma}

Lemma~\ref{lem:temporal_invariance} separates the service outcomes used to evaluate fairness from the operational decisions through which they are achieved. The two plans may differ in their repair schedules, generation and storage decisions, and timing of unmet demand. These differences can affect operational feasibility and restoration cost even when the fairness outcomes coincide.

To examine how physical access to power restricts these outcomes, let $\tau_i$ denote the earliest period in which load bus $i$ can receive positive service under any feasible restoration plan in $\Omega$. Set $\tau_i=|\mathcal{T}|+1$ if the bus cannot receive positive service during the horizon. Thus, $\gamma_{i,t}=0$ for all $t<\tau_i$. For an initially disconnected bus without available local supply, $\tau_i$ depends on the earliest feasible restoration of a supply path. The next result bounds attainable fairness independently of the operator's fairness preferences.

\begin{proposition}[Connectivity bounds on attainable fairness]
\label{prop:connectivity_bound}
For each bus $i\in\mathcal{B}^L$, we have
\[
y_i\leq\frac{|\mathcal{T}|-\tau_i+1}{|\mathcal{T}|}.
\]
Consequently,
\[
\max_{\boldsymbol{\gamma}\in\Gamma_0}\min_{i\in\mathcal{B}^L}y_i
\leq
\min_{i\in\mathcal{B}^L}
\frac{|\mathcal{T}|-\tau_i+1}{|\mathcal{T}|},
\]
and
\[
\min_{\boldsymbol{\gamma}\in\Gamma_0}
\max_{i\in\mathcal{B}^L}
\left\{
SVI_i\sum_{t\in\mathcal{T}}(1-\gamma_{i,t})
\right\}
\geq
\max_{i\in\mathcal{B}^L}
\left\{SVI_i(\tau_i-1)\right\}.
\]
\end{proposition}

Proposition~\ref{prop:connectivity_bound} identifies a physical limit on both dimensions of fairness. A bus that cannot receive power before $\tau_i$ necessarily accumulates $\tau_i-1$ periods without service. Relaxing this bound requires earlier feasible access to power, for example through additional repair capacity or local supply. Once power becomes available, however, fairness also depends on how scarce energy is distributed across buses. This raises a further question: when can the available energy support both the highest attainable minimum service level and the best vulnerability-aware outcome?

We examine this question using an aggregate allocation benchmark that retains differences in demand, vulnerability, and attainable service. Suppose demand at each bus is constant over time, and let $D_i=\sum_{t\in\mathcal{T}}p^D_{i,t}\Delta t>0$ denote the energy required to fully serve bus $i$ over the horizon. Its served energy is then $D_i y_i$. Let $\overline{E}>0$ denote the total available energy and $\overline{y}_i\in(0,1]$ a ceiling on average service. For this benchmark, assume that the attainable service vectors are exactly
\begin{equation}
\mathcal{Y}^{A}
=
\left\{
\mathbf{y}:\ 0\leq y_i\leq\overline{y}_i\ \forall i\in\mathcal{B}^L,
\quad
\sum_{i\in\mathcal{B}^L}D_i y_i\leq\overline{E}
\right\}.
\label{eq:aggregate_service_set}
\end{equation}
The service ceilings can represent restrictions such as delayed access to power. The benchmark abstracts from additional coupling caused by network flows, storage operation, and the timing of energy availability. It therefore isolates the interaction between aggregate scarcity and individual service restrictions, rather than replacing the feasible restoration set $\Omega$.

Let $m^{A}=\max_{\mathbf{y}\in\mathcal{Y}^{A}}\min_i y_i$ denote the MMF service level in this benchmark, and let $V^{A}=\min_{\mathbf{y}\in\mathcal{Y}^{A}}\allowbreak\max_i\{SVI_i|\mathcal{T}|(1-y_i)\}$ denote its vertical-fairness optimum. These benchmarks are computed over $\mathcal{Y}^{A}$ and are distinct from those computed over the full restoration model. Write $(z)_+=\max\{z,0\}$.

\begin{proposition}[Compatibility of horizontal and vertical fairness]
\label{prop:fairness_compatibility}
Suppose $SVI_i>0$ for every $i\in\mathcal{B}^L$. For the aggregate allocation benchmark in \eqref{eq:aggregate_service_set},
\begin{equation}
m^{A}
=
\min\left\{
\min_{i\in\mathcal{B}^L}\overline{y}_i,
\frac{\overline{E}}{\sum_{i\in\mathcal{B}^L}D_i}
\right\}.
\label{eq:aggregate_mmf}
\end{equation}
The value $V^{A}$ is the smallest $z\geq0$ satisfying
\begin{equation}
z\geq |\mathcal{T}|\max_{i\in\mathcal{B}^L}
\{SVI_i(1-\overline{y}_i)\},
\qquad
\sum_{i\in\mathcal{B}^L}D_i
\left(1-\frac{z}{|\mathcal{T}|SVI_i}\right)_+
\leq\overline{E}.
\label{eq:aggregate_vf}
\end{equation}
A feasible allocation attains both benchmarks if and only if
\begin{equation}
\sum_{i\in\mathcal{B}^L}D_i
\max\left\{
m^{A},\ 1-\frac{V^{A}}{|\mathcal{T}|SVI_i}
\right\}
\leq\overline{E}.
\label{eq:fairness_compatibility}
\end{equation}
In particular, if $\overline{y}_i=1$ for all buses and $0<\overline{E}<\sum_iD_i$, simultaneous attainment is possible if and only if all SVI values are equal.
\end{proposition}

Proposition~\ref{prop:fairness_compatibility} expresses compatibility as an energy requirement. Attaining MMF requires every bus to receive at least $m^{A}$, while attaining the vertical-fairness benchmark requires bus $i$ to receive at least $1-V^{A}/(|\mathcal{T}|SVI_i)$. Their maximum is the minimum service needed to meet both requirements. When the associated energy requirement exceeds availability, the two benchmarks cannot be attained together. Without individual service restrictions, scarcity and heterogeneous vulnerability necessarily create such a conflict. With service ceilings, however, a physical bottleneck can limit the benchmarks enough to permit simultaneous attainment.

An illustrative case shows that expanding physical service possibilities does not necessarily preserve this compatibility. Consider two load buses $i$ and $j$ with $D_i=D_j=1$, $\overline{E}=1$, and vulnerability scores $SVI_i=2$ and $SVI_j=1$. If $\overline{y}_i=2/5$ and $\overline{y}_j=1$, then $m^{A}=2/5$ and $V^{A}=6|\mathcal{T}|/5$. The allocation $(y_i,y_j)=(2/5,3/5)$ attains both benchmarks. Raising $\overline{y}_i$ to $1$ improves the MMF service level to $1/2$ and reduces the vertical-fairness optimum to $2|\mathcal{T}|/3$. However, their respective optimal allocations are now $(1/2,1/2)$ and $(2/3,1/3)$, so simultaneous attainment is no longer possible. The earlier allocation remains feasible, but the improved benchmarks require different distributions of service.

The same example also distinguishes an unavoidable trade-off from the selection among multiple optimal solutions. Before the service ceiling is raised, $(2/5,0)$ is vertically optimal as well, although its minimum service level is zero. Thus, a vertically optimal plan with poor horizontal fairness does not by itself establish a conflict between the two objectives. In the full restoration model, joint attainment must be checked against the operational constraints, but this distinction remains relevant when interpreting solutions with identical VFD and different DMMF values.

We next return to the full formulation to examine the role of discrete restoration decisions. Fix $\beta$ and the horizontal reference vector $\mathbf{y}^{\alpha}$ while varying the cost limit $\lambda Z^U$, where $Z^U>0$. Let $\mathcal{K}$ index the finite set of feasible combinations of line-availability, repair-crew, and battery-operating-mode decisions, and let $\Omega_k\subseteq\Omega$ contain the restoration plans with discrete decisions fixed to combination $k$. For the power-flow relaxation used in Section~\ref{sec:FeasibleRestorationSet}, fixing these decisions leaves a convex feasible set in the continuous variables. Let $\mathbf{x}$ collect the restoration decisions and let $Z(\mathbf{x})$ denote their total cost in \eqref{eq:ut_obj}. Define
\begin{equation}
\Phi_k(\lambda)
=
\min_{\mathbf{x}\in\Omega_k}
\left\{
V(\mathbf{x}):\
y_i(\mathbf{x})\geq\beta y_i^{\alpha}\ \forall i\in\mathcal{B}^L,
\quad
Z(\mathbf{x})\leq\lambda Z^U
\right\},
\label{eq:fixed_schedule_value}
\end{equation}
with $\Phi_k(\lambda)=+\infty$ when the constraints are infeasible. Let $\Phi(\lambda)$ denote the optimal value when the discrete decisions are also optimized.

\begin{proposition}[Cost limits and discrete restoration decisions]
\label{prop:restoration_value_structure}
For fixed $\beta$ and $\mathbf{y}^{\alpha}$, each $\Phi_k(\lambda)$ is nonincreasing and convex on its feasible domain. Moreover,
\[
\Phi(\lambda)=\min_{k\in\mathcal{K}}\Phi_k(\lambda).
\]
Consider three cost factors $\lambda_1<\lambda_2<\lambda_3$ for which $\Phi$ is finite, and set $\theta=(\lambda_3-\lambda_2)/(\lambda_3-\lambda_1)$. If
\begin{equation}
\Phi(\lambda_2)
>
\theta\Phi(\lambda_1)+(1-\theta)\Phi(\lambda_3),
\label{eq:schedule_change_condition}
\end{equation}
then no combination of discrete restoration decisions is optimal at both $\lambda_1$ and $\lambda_3$.
\end{proposition}

Proposition~\ref{prop:restoration_value_structure} distinguishes the effects of continuous operating flexibility from changes in discrete line-availability, repair-crew, and battery-operating-mode decisions. The condition does not, by itself, identify which of these discrete decisions must change. Under a fixed discrete restoration plan, relaxing the cost limit produces diminishing improvements in the optimal vertical-fairness objective. In particular, an interval over which the optimal value remains constant cannot be followed by a further improvement under that same plan. When discrete decisions are optimized, the model can move between different operating possibilities, and the resulting value function does not necessarily inherit the convexity of the individual functions $\Phi_k$.

A useful implication is that $\Phi(\lambda_1)=\Phi(\lambda_2)>\Phi(\lambda_3)$ satisfies \eqref{eq:schedule_change_condition}. Such a plateau followed by an improvement therefore requires different optimal discrete decisions at the two endpoint cost limits. This provides a criterion for identifying when attaining the new fairness optimum requires changing the discrete decisions of an initially optimal restoration plan. The result concerns the optimized vertical objective, and hence VFD when its benchmark and normalization remain fixed. It does not impose the same pattern on secondary measures such as DMMF or JI. Its numerical application also requires a fixed horizontal reference vector and differences established beyond rounding and solver tolerances.

\section{Numerical Results and Discussions}\label{sec:CaseStudy}
\noindent
In this section, we apply the proposed framework to the IEEE 33-bus power distribution test system to examine how efficiency, horizontal fairness, and vertical fairness shape post-disaster restoration decisions and service outcomes. We first describe the test system and experimental setup. We then use the \utilitarian{} benchmark to identify fairness concerns that can arise under efficiency-oriented restoration, followed by an analysis of the \horizontal{} and \horizontalVertical{} models and their trade-offs in terms of POF, DMMF, VFD, and the fairness indices defined in Section~\ref{sec:Formulations}. We next examine how limits on allowable restoration cost affect these trade-offs and conclude the section by discussing the practical implications of the results. The comparative optimization results use the network approximation.

The experiments use Python 3.13.18 and \texttt{Gurobi 13.0.3} \citep{gurobi} on an Apple M4 Pro with four solver threads and a 1,200-second limit per solve. The \utilitarian{} model minimizes cost directly. For HF, MMF, and HVF, we first optimize the fairness objective and then minimize cost while preserving that objective within an absolute tolerance of $10^{-6}$. All reported primary and secondary solves terminated with optimal status under the solver tolerances. The largest absolute differences between the primary incumbent and bound are $9.0\times10^{-9}$ in the main instance and $9.1\times10^{-7}$ in the identical-demand comparison. The latter allows a relative gap of $10^{-4}$ or an absolute gap of $10^{-6}$. Secondary cost solves allow an absolute gap of $10^{-6}$, equivalent to \$0.10, and achieve a maximum of $9.1\times10^{-7}$. \ref{app:solver_certificates} reports the objectives, bounds, gaps, and runtimes for both stages. Feasibility and integrality tolerances are $10^{-7}$, and reported solutions pass independent constraint checks at $10^{-6}$. Tables and figures use the same saved solutions and full-precision benchmarks, with rounding applied only for presentation.

\subsection{Test System and Experimental Setup} \label{sec:CaseStudy_Setup}
\begin{figure}[htbp!]
    \centering
    \includegraphics[width=0.8\textwidth]{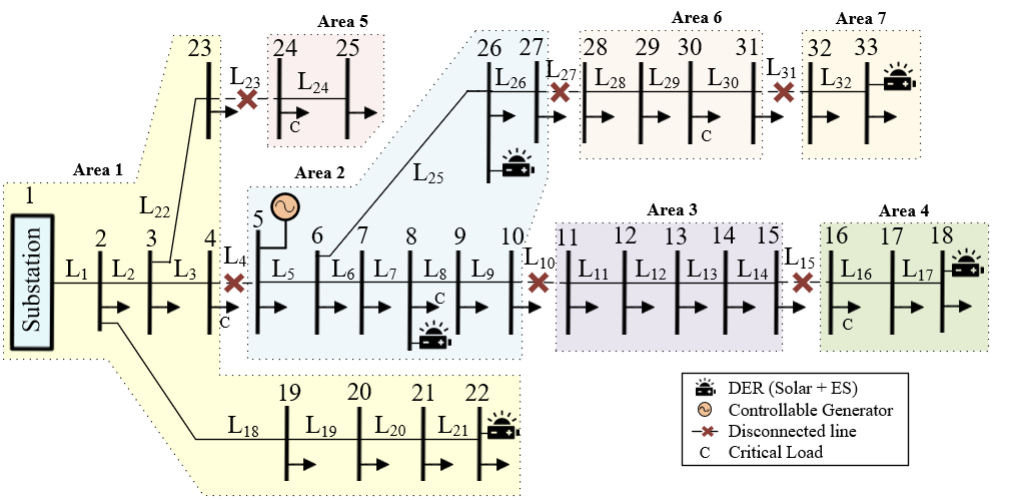}
    \caption{IEEE 33-bus test power distribution system}
    \label{fig:IEEE33}
\end{figure}

\noindent
The topology of the test system is illustrated in Figure~\ref{fig:IEEE33}. The power distribution system is connected to the transmission grid through a substation at Bus 1 and includes one controllable distributed generator at Bus 5 and five solar-plus-storage systems located at buses 8, 18, 22, 26, and 33. The system data, including 24-hour load profiles, power line characteristics, distributed energy resource parameters, diurnal solar-generation capacity, and base power and voltage values, are obtained from \cite{IEEE33Dataset}.

We set $\Delta t=1$ hour and $\mathcal{T}=\{1,\ldots,24\}$, corresponding to a 24-hour restoration horizon beginning at 8:00 AM immediately after the assumed disaster. Thus, no repair effort has been accumulated at $t=0$. The failed-line set is $\mathcal{E}^f=\{4,10,15,23,27,31\}$, where each line index corresponds to the bus pair shown in Figure~\ref{fig:IEEE33}. For example, $\ell=4$ corresponds to $(i,j)=(4,5)$. The repair requirement $C_{\ell}$ for each failed line is reported in Table~\ref{tab:repair_effort}. We assume that $\overline{C}=10$ repair crew members are available in each period and require integer crew assignments. A failed line can return to service from the period following completion of its required repair effort $C_{\ell}$. The initial line-status parameter is set to $a_{\ell,0}=0$ for $\ell\in\mathcal{E}^f$ and $a_{\ell,0}=1$ for $\ell\in\mathcal{E}\setminus\mathcal{E}^f$. As specified by the repair constraints, crew members assigned to an unfinished line remain on that line until its repair is completed.

The generation cost and VOLL are set to 73~\text{\$/MWh} and 2000~\text{\$/MWh}, respectively. Initial battery energy is zero. We further assume that power supply from the transmission system remains available throughout the restoration horizon and is sufficient to support downstream demand subject to the network and restoration constraints.

\begin{table}[htbp]
\singlespacing\small\setlength{\tabcolsep}{3pt}
    \centering
    \small
    \caption{Repair effort required to restore failed power lines, measured in repair crew-hours.}
    \label{tab:repair_effort}
    \begin{tabular}{cc}
        \hline
        \textbf{Line Index} & \textbf{Repair effort} \\ \hline
        4  & 50 \\ \hline
        10 & 20 \\ \hline
        15 & 20 \\ \hline
        23 & 35 \\ \hline
        27 & 40 \\ \hline
        31 & 40 \\ \hline
    \end{tabular}
\end{table}

The six failed lines separate the distribution system into seven areas, as shown in Figure~\ref{fig:IEEE33}. Area 1 remains connected to the transmission system and therefore retains access to grid power throughout the restoration horizon. Areas 3, 5, and 6 initially experience complete blackouts because they are disconnected from both the transmission system and local generation resources. In contrast, Areas 4 and 7 have sufficient local generation to remain self-sustaining throughout the horizon, while Area 2 can maintain at least partial service temporarily using its distributed generation and solar-plus-storage resources. Across the seven areas, aggregate power demand is ranked from largest to smallest as Areas 1, 2, 5, 3, 6, 4, and 7, respectively.

\subsection{Utilitarian-Based Restoration} \label{sec:CaseStudy_Efficiency}
\begin{figure}[htbp!]
    \captionsetup[subfigure]{font=footnotesize, skip=1pt}
    \centering

    \begin{subfigure}[t]{0.32\textwidth}
        \centering
        \includegraphics[width=\textwidth]{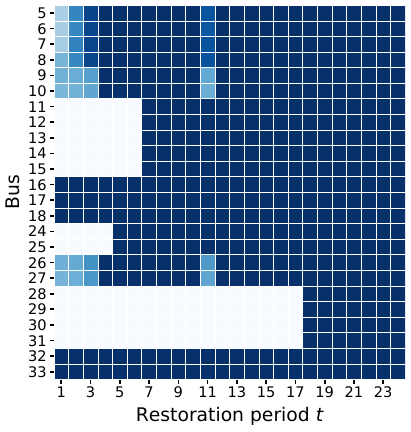}
        \caption{Utilitarian Benchmark}
        \label{fig:panel_HF_a}
    \end{subfigure}\hfill
    \begin{subfigure}[t]{0.32\textwidth}
        \centering

    \end{subfigure}\hfill
    \begin{subfigure}[t]{0.32\textwidth}
        \centering
       \includegraphics[width=\linewidth,height=3cm,keepaspectratio]{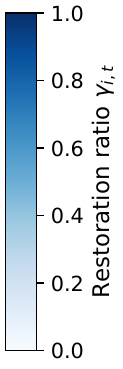}
    \end{subfigure}

    \vspace{0.1em}

    \begin{subfigure}[t]{0.32\textwidth}
        \centering
        \includegraphics[width=\textwidth]{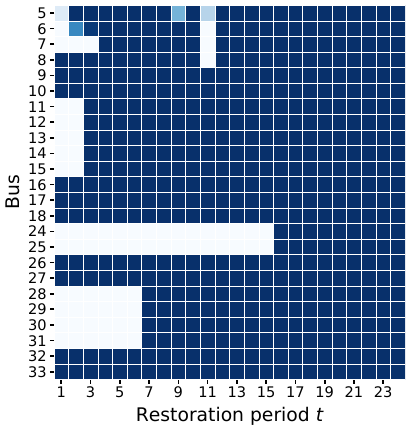}
        \caption{Horizontal Fairness, $\alpha=0.25$}
        \label{fig:panel_HF_b}
    \end{subfigure}\hfill
    \begin{subfigure}[t]{0.32\textwidth}
        \centering
        \includegraphics[width=\textwidth]{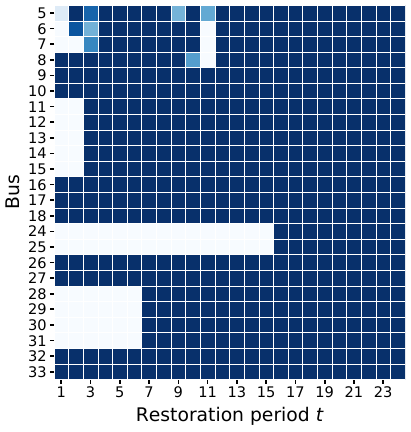}
        \caption{Horizontal Fairness, $\alpha=0.5$}
        \label{fig:panel_HF_c}
    \end{subfigure}\hfill
    \begin{subfigure}[t]{0.32\textwidth}
        \centering
        \includegraphics[width=\textwidth]{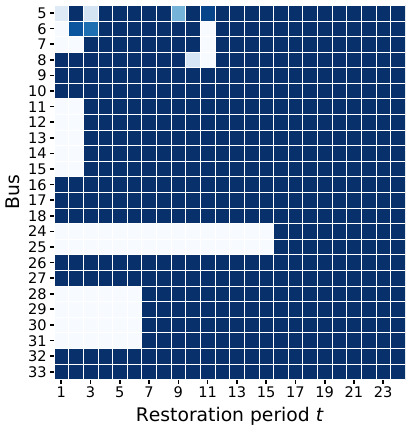}
        \caption{Proportional Fairness, $\alpha=1$}
        \label{fig:panel_HF_d}
    \end{subfigure}

    \begin{subfigure}[t]{0.32\textwidth}
        \centering
        \includegraphics[width=\textwidth]{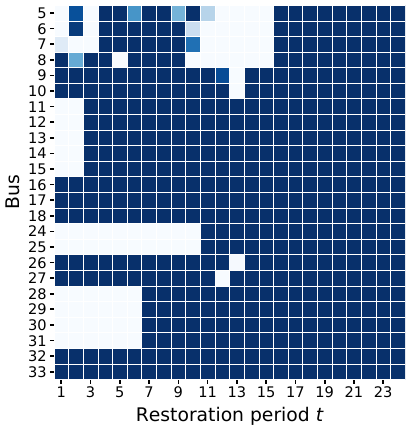}
        \caption{Horizontal Fairness, $\alpha=2$}
        \label{fig:panel_HF_e}
    \end{subfigure}\hfill
    \begin{subfigure}[t]{0.32\textwidth}
        \centering
        \includegraphics[width=\textwidth]{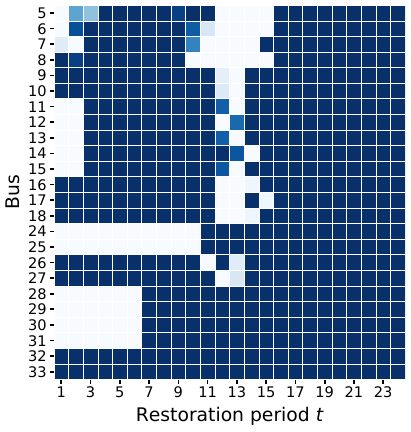}
        \caption{Horizontal Fairness, $\alpha = 5$}
        \label{fig:panel_HF_f}
    \end{subfigure}\hfill
    \begin{subfigure}[t]{0.32\textwidth}
        \centering
        \includegraphics[width=\textwidth]{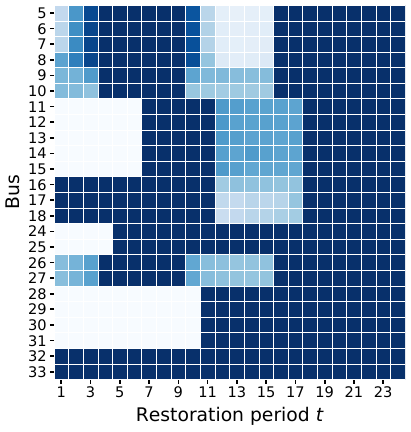}
        \caption{Max-Min Fairness}
        \label{fig:panel_HF_g}
    \end{subfigure}

    \caption{Bus-level restoration ratios $\gamma_{i,t}$ over the 24-hour restoration horizon under the \utilitarian{} benchmark and \horizontal{} models with different values of $\alpha$. Buses in Area 1 are not shown.}
    \label{fig:gamma_compare_HE}
\end{figure}

\begin{figure}[ht!]
    \centering
    \includegraphics[width=\textwidth]{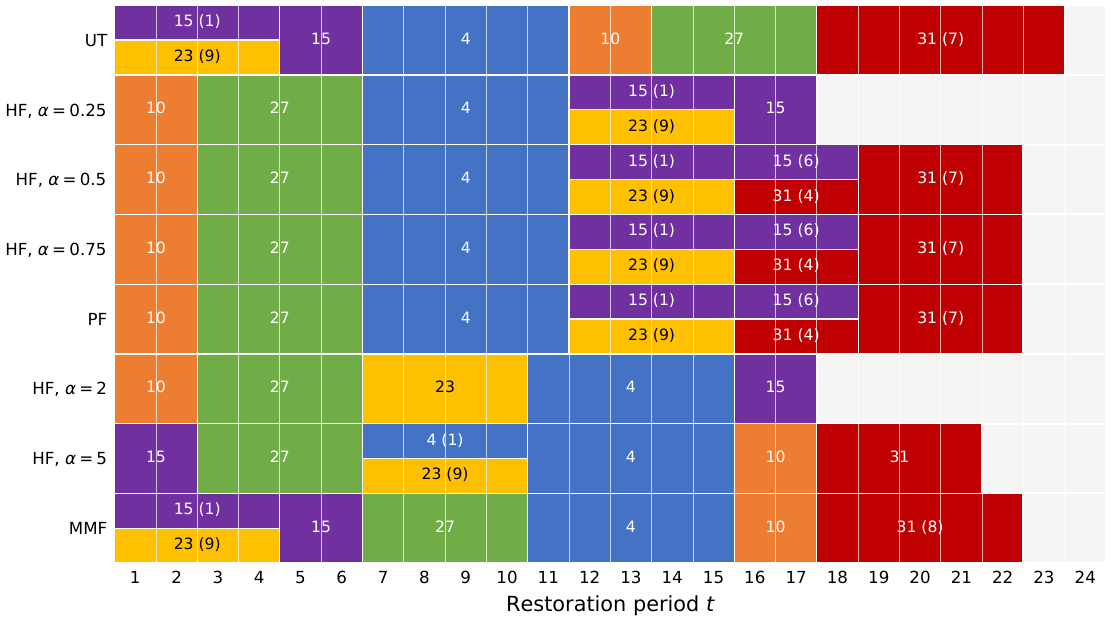}
    \caption{Crew assignment over the 24-hour restoration horizon under the \utilitarian{} benchmark and \horizontal{} models. Colors and numbers identify the assigned lines. Parentheses indicate crew counts below ten, and multiple entries indicate simultaneous repairs. The color version of the figure is available online.}
    \label{fig:crew_compare_UT_HE}
\end{figure}

\noindent
We first examine the restoration plan obtained from the \utilitarian{} benchmark to identify fairness concerns that can arise when restoration decisions are driven by system-wide efficiency. The bus-level service outcomes over the 24-hour horizon are shown in Figure~\ref{fig:panel_HF_a}, where darker cells indicate higher fractions of demand served and lighter cells indicate greater unmet demand. The corresponding repair-crew schedule is shown in the first row of Figure~\ref{fig:crew_compare_UT_HE}. Crews may be divided between failed lines. Accordingly, Figure~\ref{fig:crew_compare_UT_HE} identifies the assigned lines and indicates crew counts when fewer than ten members work on a line.

Under the \utilitarian{} benchmark, Lines 23 and 15 are repaired first, reconnecting Area 5 and connecting Area 3 to Area 4's local energy resources, respectively. These two areas - Area 5 and 3 - have the third- and fourth-largest aggregate demand among the seven areas. The repair sequence then shifts toward Area 2, which has local distributed generation and solar-plus-storage resources but still experiences substantial unmet demand because of its relatively high demand. Area 6, which has lower aggregate demand and no local generation resources, is consequently restored later and remains in complete blackout for 17 hours. Thus, while the \utilitarian{} solution minimizes total restoration cost, its repair sequence can produce substantial differences in service outcomes across areas. In particular, prioritizing repairs with greater system-wide benefits can delay the restoration of lower-demand areas that depend entirely on reconnection for power supply.
This outcome motivates the \horizontal{} models examined in the next subsection.

Figure~\ref{fig:panel_HF_a} also shows periods of partial load shedding across several groups of buses beginning around $t=10$. This pattern reflects the time-varying availability of solar generation and the limited capacity of the associated storage systems. Solar generation peaks around 1:00 PM and declines substantially after 5:00 PM. Period $t=1$ begins at 8:00 AM, so period $t=10$ covers 5:00--6:00 PM, when solar generation is already declining. Although the storage systems can shift some surplus solar energy to later periods, their stored energy is not sufficient to eliminate evening load shedding in all affected areas. Reconnection to the transmission system therefore remains important for maintaining service after local solar and stored energy become insufficient.

\subsection{Horizontal Fairness and Efficiency Trade-offs}
\label{sec:CaseStudy_HE}
\noindent
This section examines how the \horizontal{} models change restoration decisions relative to the \utilitarian{} benchmark. We solve the \horizontal{} formulation for different values of $\alpha \in \{0.25,\allowbreak 0.5,\allowbreak 0.75,\allowbreak 1,\allowbreak 2,\allowbreak 5,\allowbreak\infty\}$, including PF $(\alpha=1)$ and MMF $(\alpha=\infty)$. Table~\ref{tab:price_compare_UT_HE} summarizes the resulting total cost, cost increase, POF, DMMF, and JI. Figures~\ref{fig:panel_HF_b}--\ref{fig:panel_HF_g} and~\ref{fig:crew_compare_UT_HE} show the corresponding bus-level service patterns and repair-crew schedules.

\begin{table}[htbp]
\singlespacing\small\setlength{\tabcolsep}{3pt}
    \centering
    \caption{Comparison between models with unlimited budget ($\lambda=\infty$)}
    \label{tab:price_compare_UT_HE}
    \begin{tabular*}{\linewidth}{@{\extracolsep{\fill}}lccccc@{}}
        \toprule
        \textbf{Model} & \shortstack{\textbf{Cost}\\($10^5$\$)} & \shortstack{\textbf{Cost}\\\textbf{increase (\%)}}
        & \textbf{POF} & \textbf{DMMF} & \textbf{JI}\\
        \midrule
        UT & 0.1471 & - & 0.0000 & 0.5000 & 0.9328\\
        HF, $\alpha = 0.25$ & 0.2266 & 54.0\% & 0.0596 & 0.3571 & 0.9696\\
        HF, $\alpha = 0.5$ & 0.2266 & 54.0\% & 0.0596 & 0.3571 & 0.9698\\
        HF, $\alpha = 0.75$ & 0.2266 & 54.0\% & 0.0596 & 0.3571 & 0.9699\\
        PF & 0.2266 & 54.0\% & 0.0596 & 0.3571 & 0.9699\\
        HF, $\alpha = 2$ & 0.2761 & 87.7\% & 0.0967 & 0.0000 & 0.9754\\
        HF, $\alpha = 5$ & 0.2574 & 75.0\% & 0.0827 & 0.0000 & 0.9811\\
        MMF (HF, $\alpha = \infty$) & 0.2241 & 52.4\% & 0.0577 & 0.0000 & 0.9667\\
        \bottomrule
    \end{tabular*}
\end{table}

As shown in Figure~\ref{fig:crew_compare_UT_HE}, compared with the \utilitarian{} benchmark, the MMF model moves the repair of Line 27 earlier in the horizon, immediately after Lines 23 and 15. This shift substantially improves service in Area 6, reducing its complete-blackout duration from 17 hours to 10 hours, with essentially full service thereafter. The improvement comes with a shift in the scarcity burden toward Area 2, where some buses experience greater load shedding during periods 10 through 15, as shown in Figure~\ref{fig:gamma_compare_HE}. The MMF solution therefore does not eliminate scarcity, but changes how service shortfalls are distributed across the affected areas to improve the outcomes of less-served buses. This comparison illustrates the connectivity restriction underlying Proposition~\ref{prop:connectivity_bound}. Conditional on the reported repair schedules, a 17-hour initial blackout limits average service over the 24-hour horizon to $7/24$, whereas a 10-hour blackout raises this ceiling to $14/24$. Redistributing power after reconnection cannot remove the service loss accumulated during the initial blackout. Improving this component of fairness therefore requires earlier access to power through restoration of a supply path or additional local supply.

Figure~\ref{fig:crew_compare_UT_HE} also shows that the tested cases with $\alpha\leq1$ share the same first three repair priorities, Lines 10, 27, and 4, although later crew allocations can differ. Larger values of $\alpha$ change these priorities and produce more pronounced differences in the outcomes reported in Table~\ref{tab:price_compare_UT_HE}. Figure~\ref{fig:gamma_compare_HE} further illustrates how bus-level service changes as $\alpha$ places progressively greater emphasis on balancing service outcomes across buses.

One might expect a stronger emphasis on horizontal fairness to monotonically reduce efficiency. Table~\ref{tab:price_compare_UT_HE}, however, shows that this is not always the case. In particular, after secondary cost minimization, the MMF solution has lower total cost and POF than the intermediate \horizontal{} solutions considered here. This nonmonotonic relationship arises because the \horizontal{} objectives are defined using the average fraction of demand served at each bus, whereas system efficiency depends on the aggregate amount of demand served. As a result, the fairness objectives treat proportional service outcomes across buses independently of their demand magnitudes, while POF weights these outcomes by bus-level demand.

\begin{figure}[ht!]
    \centering
    \includegraphics[width=\textwidth]{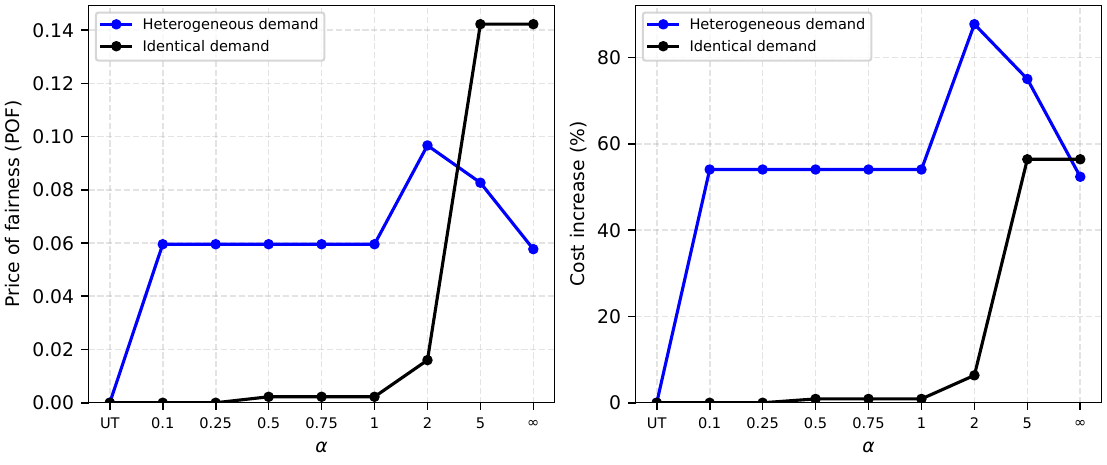}
    \caption{POF and cost increase relative to the \utilitarian{} benchmark across values of $\alpha$ under heterogeneous and identical bus-level demand}
    \label{fig:Alpha_POF_demand}
\end{figure}

To examine the role of demand heterogeneity in this result, we consider a variant in which all buses have identical demand in each period, with active and reactive demand set equal to their corresponding period-specific averages across the 32 load buses in the original dataset. Each demand scenario uses its own utilitarian benchmark, shown as UT in Figure~\ref{fig:Alpha_POF_demand}. The comparison also includes $\alpha=0.1$. As shown in Figure~\ref{fig:Alpha_POF_demand}, POF is nondecreasing with $\alpha$ for the tested values in this case, up to numerical tolerances, and levels off at 0.1423 from $\alpha=5$ to MMF. Both solutions have a cost increase of 56.4\% after secondary cost minimization. The contrast with the original instance indicates that the efficiency cost of horizontal fairness depends not only on the chosen value of $\alpha$, but also on the distribution of demand across the network. Consequently, a stronger emphasis on horizontal fairness does not necessarily translate into a proportionally larger efficiency loss when bus-level demands are heterogeneous.

In the original instance, DMMF exhibits a different pattern. As $\alpha$ increases, DMMF decreases and reaches zero at $\alpha=2$, indicating that the minimum average service level under the $\alpha=2$ and $\alpha=5$ solutions matches that of the MMF benchmark. These solutions can nevertheless differ from MMF in their repair schedules and in the distribution of service among buses above the minimum service level. The differences in JI reported in Table~\ref{tab:price_compare_UT_HE} capture these broader changes in the distribution of bus-level service.

\subsection{Vulnerability-Aware Restoration}
\label{sec:CaseStudy_VE}
\noindent
We next examine how vulnerability-aware restoration changes the service-allocation patterns obtained from the \horizontal{} models. While horizontal fairness seeks to balance service outcomes across buses, it does not distinguish among buses with different vulnerability levels. The \horizontalVertical{} model accounts for these differences through SVI-weighted unmet demand, while $\beta$ controls how closely the service levels of the corresponding \horizontal{} solution must be preserved. Lower values of $\beta$ provide greater flexibility to prioritize vertical fairness, whereas $\beta=1$ requires the service level of each bus to be at least that obtained from the corresponding \horizontal{} solution.

\begin{table}[htbp]
\singlespacing\small\setlength{\tabcolsep}{3pt}
    \centering
    \caption{Synthetic SVI value per bus in the IEEE 33-bus test power distribution system.}
    \label{tab:bus_svi_mapping}
    \begin{tabular*}{\linewidth}{@{\extracolsep{\fill}}cccccccc@{}}
    \toprule
    \textbf{Bus ID} & \textbf{SVI} & \hspace{0.3cm} & \textbf{Bus ID} & \textbf{SVI} & \hspace{0.3cm} & \textbf{Bus ID} & \textbf{SVI} \\
    \midrule
    1  & 5 & & 12 & 2 & & 23 & 4 \\
    2  & 2 & & 13 & 3 & & 24 & 3 \\
    3  & 5 & & 14 & 4 & & 25 & 2 \\
    4  & 5 & & 15 & 4 & & 26 & 5 \\
    5  & 1 & & 16 & 3 & & 27 & 1 \\
    6  & 4 & & 17 & 3 & & 28 & 5 \\
    7  & 1 & & 18 & 5 & & 29 & 1 \\
    8  & 5 & & 19 & 1 & & 30 & 1 \\
    9  & 1 & & 20 & 3 & & 31 & 4 \\
    10 & 4 & & 21 & 3 & & 32 & 4 \\
    11 & 1 & & 22 & 5 & & 33 & 1 \\
    \bottomrule
    \end{tabular*}
\end{table}

\begin{figure}[htbp!]
    \captionsetup[subfigure]{font=footnotesize, skip=1pt}
    \centering

    \begin{subfigure}[t]{0.33\textwidth}
        \centering
        \includegraphics[width=\textwidth]{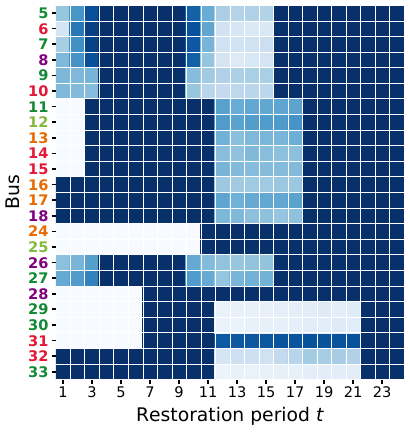}
        \caption{Horizontal-vertical Fairness,\\ $\beta=0$}
        \label{fig:panel_HVF_a}
    \end{subfigure}\hfill
    \begin{subfigure}[t]{0.33\textwidth}
        \centering
        \includegraphics[width=\textwidth]{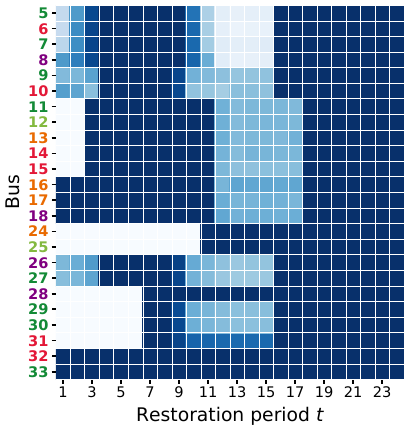}
        \caption{Horizontal-vertical Fairness,\\ $\alpha=2, \>\beta=0.75$}
        \label{fig:panel_HVF_b}
    \end{subfigure}\hfill
    \begin{subfigure}[t]{0.33\textwidth}
        \centering
       \includegraphics[width=\textwidth]{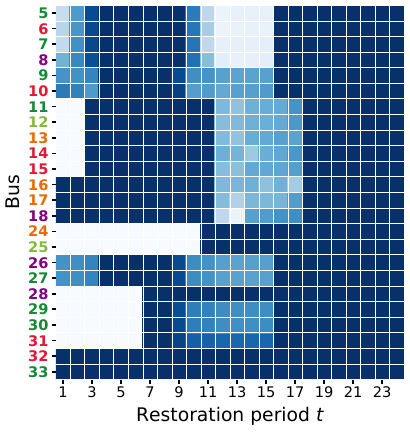}
       \caption{Horizontal-vertical Fairness,\\ $\alpha=2,\>\beta=0.85$}
        \label{fig:panel_HVF_c}
    \end{subfigure}

    \vspace{0.1em}

    \begin{subfigure}[t]{0.33\textwidth}
        \centering
        \includegraphics[width=\textwidth]{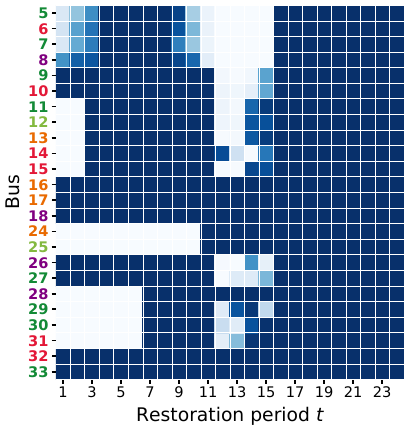}
        \caption{Horizontal-vertical Fairness,\\ $\alpha=2,\>\beta=0.9$}
        \label{fig:panel_HVF_d}
    \end{subfigure}\hfill
    \begin{subfigure}[t]{0.33\textwidth}
        \centering
        \includegraphics[width=\textwidth]{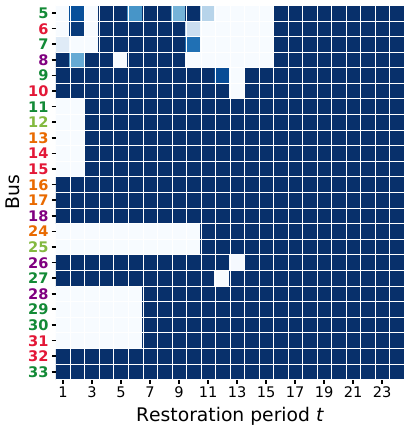}
        \caption{Horizontal Fairness, $\alpha=2$}
        \label{fig:panel_HVF_e}
    \end{subfigure}\hfill
    \begin{subfigure}[t]{0.33\textwidth}
        \centering
        \includegraphics[width=\linewidth,height=4cm,keepaspectratio]{Restoration_colorbar.pdf}

    \end{subfigure}

    \caption{Comparison of the ratio of power supplied ($\gamma_{i,t}$) to each bus for \ref{fig:panel_HVF_a}: sole focus on vertical fairness, \ref{fig:panel_HVF_b} to \ref{fig:panel_HVF_d}: a range of $\beta$ values where $\beta \neq 1$ and $\alpha = 2$, and \ref{fig:panel_HVF_e}: the horizontal reference corresponding to $\beta=1$. Color on the y-axis ticks indicates the SVI of each bus, with purple to green matching highest to lowest SVI. Buses in Area 1 are not shown. The color version of the figure is available online.}
    \label{fig:gamma_compare_HVE}
\end{figure}

The SVI scores assigned to the 33 buses are randomly generated integers from 1 to 5 and are reported in Table~\ref{tab:bus_svi_mapping}. Figure~\ref{fig:gamma_compare_HVE} compares the bus-level restoration ratios under selected \horizontalVertical{} models with those under the \horizontal{} model with $\alpha=2$. The results illustrate that vulnerability-aware restoration can change power allocation even when the repair schedule remains unchanged. This effect is particularly visible within connected areas with limited available power. For example, during periods 10 through 15 in plots from Figures~\ref{fig:panel_HVF_b} and \ref{fig:panel_HVF_c}, Areas 2 and 6 share the available solar-plus-storage resources. As $\beta$ decreases, the model gives greater priority to buses with higher SVI values. Bus 28, with $SVI=5$, remains fully served during these periods, while lower-SVI buses such as Buses 29 and 30, both with $SVI=1$, absorb a larger share of the load shedding. Because plots from Figures~\ref{fig:panel_HVF_b} and \ref{fig:panel_HVF_c} have the same repair-crew schedule, their differences isolate the effect of $\beta$ on power allocation within the same connected portion of the network. For $\beta=0.9$, Figure~\ref{fig:panel_HVF_d} gives Bus 31 a higher average service level than the lower-SVI Buses 29 and 30, although their period-by-period allocations vary. These results show that the \horizontalVertical{} model can improve vertical fairness through both repair priorities and the allocation of scarce power among buses that are already connected.

\begin{figure}[ht!]
    \centering
    \includegraphics[width=\textwidth]{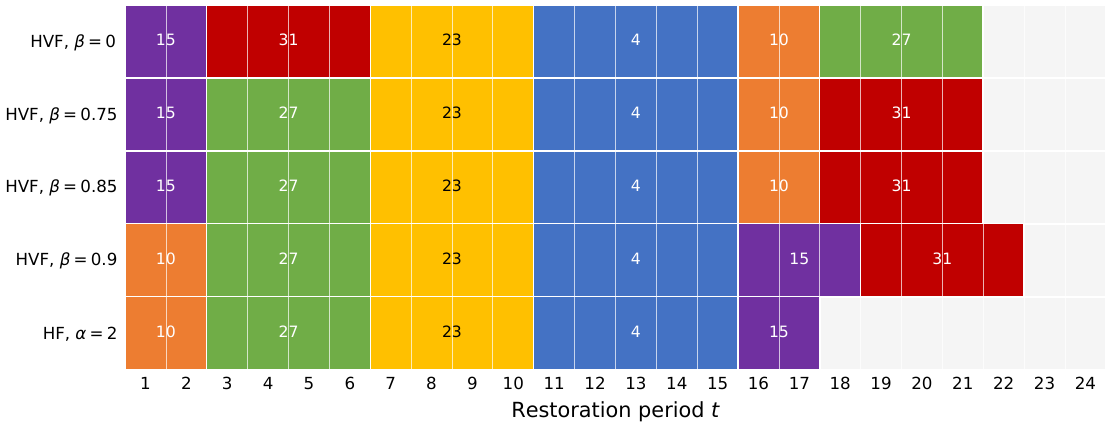}
    \caption{Crew assignment over the 24-hour restoration horizon under the \horizontalVertical{} models. Colors and numbers identify the assigned lines. Parentheses indicate crew counts below ten, and multiple entries indicate simultaneous repairs. The color version of the figure is available online.}
    \label{fig:crew_compare_UT_HVE}
\end{figure}

Table~\ref{tab:price_compare_HVE} summarizes the corresponding efficiency and fairness trade-offs. Higher restoration cost does not necessarily yield better vertical fairness. The \horizontal{} model with $\alpha=2$ has the highest cost among the reported cases while retaining a positive VFD. In contrast, the \horizontalVertical{} models with $\beta\in\{0.75,0.85,0.9\}$ attain VFD equal to zero while maintaining DMMF equal to zero. Thus, these solutions attain the vertical-fairness benchmark without reducing the minimum service level relative to the MMF benchmark. Moreover, the solutions with $\beta=0.75 \text{ to } 0.9$ achieve these outcomes at a lower cost than the \horizontal{} model with $\alpha=2$. For this range of $\beta$, the model can therefore improve vulnerability-aware allocation while preserving the main horizontal-fairness outcome obtained in Section~\ref{sec:CaseStudy_HE}.

\begin{table}[htbp]
\singlespacing\small\setlength{\tabcolsep}{3pt}
    \centering
    \caption{Comparison between models with unlimited budget}
    \label{tab:price_compare_HVE}
    \begin{tabular*}{\linewidth}{@{\extracolsep{\fill}}lccccc@{}}
        \toprule
        \textbf{Model} & \shortstack{\textbf{Cost}\\($10^5$\$)} & \shortstack{\textbf{Cost}\\\textbf{increase (\%)}} & \textbf{POF} &\textbf{DMMF} &\textbf{VFD}\\
        \midrule
        HVF ($\beta = 0$) & 0.2594 & 76.4\% & 0.0842 & 0.3634 & 0.0000\\
        HVF ($\alpha = 2, \beta = 0.75$) & 0.2597 & 76.5\% & 0.0844 & 0.0000 & 0.0000\\
        HVF ($\alpha = 2, \beta = 0.85$) & 0.2597 & 76.5\% & 0.0844 & 0.0000 & 0.0000\\
        HVF ($\alpha = 2, \beta = 0.9$) & 0.2735 & 85.9\% & 0.0947 & 0.0000 & 0.0000\\
        HF $(\alpha = 2)$ & 0.2761 & 87.7\% & 0.0967 & 0.0000 & 0.0826\\
        \bottomrule
    \end{tabular*}
\end{table}

The distinction between compatibility of fairness benchmarks and
selection among optimal allocations, illustrated by
Proposition~\ref{prop:fairness_compatibility}, is relevant to these
results. The solutions with $\beta\in\{0.75,0.85,0.9\}$ attain both
DMMF and VFD equal to zero. These solutions remain feasible when
$\beta=0$, since removing the horizontal-service requirements enlarges
the feasible set, and they retain the optimal vertical objective.
Consequently, the DMMF of 0.3634 reported for the cost-minimizing
vertically optimal solution at $\beta=0$ does not represent an unavoidable sacrifice for
attaining vertical fairness. It shows that the vertical objective and secondary cost minimization
do not enforce the MMF service level. In this instance, the horizontal-service requirements
select vertically optimal plans that also preserve the MMF service
level.

These results also clarify the mechanisms through which vertical fairness can be improved. Because restoration cost is strongly influenced by discrete repair decisions and the resulting grid connectivity, different models can produce similar cost and POF values. Yet, for a given repair schedule and pattern of available power, the operator may still have flexibility in allocating power across connected buses. The \horizontalVertical{} model uses this flexibility to prioritize higher-SVI buses during periods of scarcity. Vulnerability-aware restoration can therefore improve vertical fairness not only through changes in repair priorities, but also through the allocation of available power within connected areas.

\begin{figure}[ht!]
    
    \centering
    \includegraphics[width=0.55\textwidth]{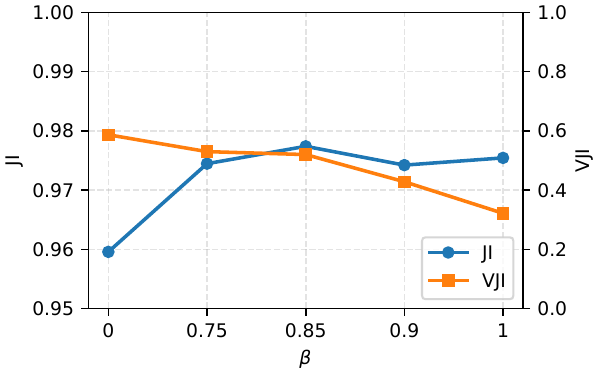}
    \caption{The Jain Index (JI) and the Vulnerability-adjusted Jain Index (VJI) for the \horizontalVertical{} models with $\alpha=2$. The horizontal reference is shown at $\beta=1$.
    }
    \label{fig:J_MJ_HVE}
    
\end{figure}

Finally, Figure~\ref{fig:J_MJ_HVE} shows how JI and VJI change with $\beta$ for $\alpha=2$. As $\beta$ decreases, VJI generally increases, indicating a more balanced distribution of vulnerability-weighted unmet service across buses. JI remains relatively stable over a range of $\beta$ values but eventually declines as greater flexibility is given to vulnerability-aware allocation. The results therefore show that vertical-fairness improvements can initially be achieved while largely preserving the overall distribution of service across buses, although stronger emphasis on vertical fairness can eventually reduce horizontal fairness.

\subsection{Resource Scarcity and Fairness Trade-offs}
\noindent
The preceding sections evaluate \horizontal{} and \horizontalVertical{} models without imposing an explicit upper bound on total restoration cost. In practice, however, system operators may face budgetary or efficiency-tolerance limits that restrict the additional cost that can be incurred to improve fairness. We therefore examine how the proposed models perform under different values of the budget factor $\lambda$, where smaller values impose tighter limits on the allowable restoration cost relative to the \utilitarian{} benchmark.

The unconstrained results in Tables~\ref{tab:price_compare_UT_HE} and~\ref{tab:price_compare_HVE} show that fairness-aware restoration can increase total cost substantially relative to the \utilitarian{} benchmark, with increases ranging from approximately 52\% to 88\% among the reported \horizontal{} and \horizontalVertical{} solutions. Based on these results, we consider three budget-factor values, $\lambda \in \{1.25,1.4,1.5\}$, representing increasingly relaxed cost limits. Tables~\ref{tab:budget_HE} and~\ref{tab:budget_HVE} report the resulting outcomes for the \horizontal{} and \horizontalVertical{} models, respectively, together with the unconstrained case $\lambda=\infty$.

\begin{table}[htbp]
\singlespacing\small\setlength{\tabcolsep}{3pt}
    \centering
    \caption{Comparison across \horizontal{} models under different budget-factor values}
    \label{tab:budget_HE}
    \renewcommand{\arraystretch}{1.12}
    \begin{tabular*}{\linewidth}{@{\extracolsep{\fill}}llccccccc@{}}
        \toprule
        $\alpha$ & $\lambda$ & \shortstack{\textbf{Cost}\\($10^5$\$)} & \shortstack{\textbf{Cost}\\\textbf{increase (\%)}} & \textbf{POF} & \textbf{DMMF} & \textbf{VFD} & \textbf{JI} & \textbf{VJI} \\
        \midrule
        1 & 1.25 & 0.160 & 8.64\% & 0.010 & 0.286 & 0.444 & 0.956 & 0.181 \\
        1 & 1.4 & 0.160 & 8.64\% & 0.010 & 0.286 & 0.444 & 0.956 & 0.181 \\
        1 & 1.5 & 0.160 & 8.64\% & 0.010 & 0.286 & 0.444 & 0.956 & 0.181 \\
        1 & $\infty$ & 0.227 & 54.03\% & 0.060 & 0.357 & 0.167 & 0.970 & 0.240 \\
        \midrule
        2 & 1.25 & 0.160 & 8.64\% & 0.010 & 0.286 & 0.444 & 0.956 & 0.181 \\
        2 & 1.4 & 0.160 & 8.64\% & 0.010 & 0.286 & 0.444 & 0.956 & 0.181 \\
        2 & 1.5 & 0.160 & 8.64\% & 0.010 & 0.286 & 0.444 & 0.956 & 0.181 \\
        2 & $\infty$ & 0.276 & 87.71\% & 0.097 & 0.000 & 0.083 & 0.975 & 0.322 \\
        \midrule
        10 & 1.25 & 0.160 & 8.64\% & 0.010 & 0.286 & 0.444 & 0.956 & 0.181 \\
        10 & 1.4 & 0.160 & 8.64\% & 0.010 & 0.286 & 0.444 & 0.956 & 0.181 \\
        10 & 1.5 & 0.160 & 8.64\% & 0.010 & 0.286 & 0.444 & 0.956 & 0.181 \\
        10 & $\infty$ & 0.257 & 75.00\% & 0.083 & 0.000 & 0.004 & 0.981 & 0.453 \\
        \bottomrule
    \end{tabular*}
\end{table}

Table~\ref{tab:budget_HE} shows a plateau over the three tested finite cost limits. For $\alpha\in\{1,2,10\}$ and $\lambda\in\{1.25,1.4,1.5\}$, the reported solutions have total cost approximately $0.159796\times10^5$ dollars, an 8.64\% increase over the \utilitarian{} benchmark, with DMMF equal to 0.286 and VFD equal to 0.444 after rounding. These cost limits are nonbinding at the reported solutions, but they still exclude other restoration plans: the finite-budget results differ from their unconstrained counterparts. For each fixed $\alpha$, the full-precision utility values and average service vectors agree across the three finite limits. Thus, additional cost tolerance within this tested range does not improve the horizontal-fairness objectives, whereas removing the cost constraint permits different service outcomes. The plateau should not be interpreted as a general invariance to either $\alpha$ or the cost limit.

\begin{table}[htbp]
\singlespacing\small\setlength{\tabcolsep}{3pt}
    \centering
    \caption{Comparison across \horizontalVertical{} models under different budget-factor values ($\alpha = 2$)}
    \label{tab:budget_HVE}
    \renewcommand{\arraystretch}{1.12}
    \begin{tabular*}{\linewidth}{@{\extracolsep{\fill}}llccccccc@{}}
        \toprule
        $\beta$ & $\lambda$ & \shortstack{\textbf{Cost}\\($10^5$\$)} & \shortstack{\textbf{Cost}\\\textbf{increase (\%)}} & \textbf{POF} & \textbf{DMMF} & \textbf{VFD} & \textbf{JI} & \textbf{VJI} \\
        \midrule
        0 & 1.25 & 0.155 & 5.59\% & 0.006 & 0.357 & 0.389 & 0.911 & 0.277 \\
        0 & 1.4 & 0.192 & 30.39\% & 0.034 & 0.357 & 0.333 & 0.925 & 0.343 \\
        0 & 1.5 & 0.220 & 49.43\% & 0.054 & 0.357 & 0.167 & 0.947 & 0.299 \\
        0 & $\infty$ & 0.259 & 76.35\% & 0.084 & 0.363 & 0.000 & 0.960 & 0.587 \\
        \midrule
        0.75 & 1.25 & 0.160 & 8.64\% & 0.010 & 0.286 & 0.444 & 0.956 & 0.181 \\
        0.75 & 1.4 & 0.160 & 8.64\% & 0.010 & 0.286 & 0.444 & 0.956 & 0.181 \\
        0.75 & 1.5 & 0.160 & 8.64\% & 0.010 & 0.286 & 0.444 & 0.956 & 0.181 \\
        0.75 & $\infty$ & 0.260 & 76.54\% & 0.084 & 0.000 & 0.000 & 0.974 & 0.530 \\
        \bottomrule
    \end{tabular*}
\end{table}

The \horizontalVertical{} results in Table~\ref{tab:budget_HVE} distinguish the effect of the budget from that of the protected horizontal reference. For $\beta=0$, VFD decreases from 0.389 to 0.333 and then to 0.167 as $\lambda$ increases from 1.25 to 1.4 and 1.5. The unrestricted case attains the vertical-fairness benchmark. For $\beta=0.75$, the finite-budget solutions instead have VFD equal to 0.444 and DMMF equal to 0.286 after rounding, matching the corresponding $\alpha=2$ horizontal reference on these measures. Their cost is approximately 8.64\% above the \utilitarian{} benchmark. Thus, under the tested finite limits, the protected service requirements prevent the improvements in the worst vulnerability-weighted shortfall attained at $\beta=0$. In the unrestricted case, the $\beta=0.75$ solution attains both fairness benchmarks.

At each finite cost factor, the horizontal reference and the corresponding horizontal-vertical-fairness solution both have vertical objective 70, consistent with \eqref{eq:hf_reference_dominance}. With $V^{VF}=30$, their common VFD is $(70-30)/(120-30)=0.444\ldots$. Without the cost constraint, the horizontal-vertical-fairness solution improves VFD from approximately 0.083 for the horizontal reference to zero.

The $\beta=0$ results also provide a numerical application of Proposition~\ref{prop:restoration_value_structure}. The vertical objectives at $\lambda=1.25$, $1.4$, and $1.5$ are 65, 60, and 45, respectively, up to solver tolerances. Since $1.4=0.4(1.25)+0.6(1.5)$ and $60>0.4(65)+0.6(45)=53$, the reported values satisfy the proposition's inequality. The archived middle lower bound is 60, and the feasible endpoint values give a weighted upper bound of 53, leaving a margin greater than 6.99999 after numerical tolerances. The fixed VFD normalization preserves this inequality. Under the proposition's assumptions, no single combination of discrete restoration decisions can be optimal at both endpoint cost limits. This does not identify a repair-schedule change specifically. For $\beta>0$, comparisons across cost limits must additionally account for any change in the protected horizontal reference vector. In the present finite-budget cases, the reported horizontal reference and the resulting vertical objective remain unchanged.

\subsection{Practical Implications}
\noindent

The numerical results yield several practical implications beyond showing that fairness-aware restoration can improve service outcomes. The main value of the proposed framework is diagnostic: it helps system operators identify whether unfairness is driven primarily by repair sequencing, limited local power availability, vulnerability-insensitive dispatch, or tight efficiency-tolerance limits. This distinction matters because each source of unfairness calls for a different operational response.

Under the \utilitarian{} benchmark, Area 6 remains in complete blackout for 17 hours, compared with 10 hours under MMF. The improvement requires earlier restoration of a supply path, since reallocating power cannot compensate for the absence of access to a source. This comparison shows why operators should assess bus-level outage durations alongside total restoration cost when setting repair priorities. Local generation and storage can support some disconnected areas temporarily, but repair planning must also account for when those resources become insufficient.

The vulnerability-aware results show that protecting minimum service and prioritizing vulnerable communities can be compatible. For $\alpha=2$, the solutions with $\beta\in\{0.75,0.85,0.9\}$ attain both fairness benchmarks at a lower cost than the corresponding \horizontal{} reference. Moreover, the $\beta=0.75$ and $0.85$ cases share the same repair schedule but allocate power differently, illustrating the flexibility available after repair priorities have been set. In this instance, operators can improve vulnerability-aware service while retaining the optimal minimum service level. Explicit service protections help identify such allocations, since minimizing the vertical objective and then cost does not itself ensure that the MMF service level is preserved.

Additional cost tolerance does not produce uniform fairness gains. The horizontal-fairness objectives remain unchanged across the three tested finite cost limits. For $\beta=0$, however, the vertical objective improves, and Proposition~\ref{prop:restoration_value_structure} establishes that no single combination of discrete decisions is optimal at both endpoint limits. With $\beta=0.75$, the protected service requirements leave the vertical objective unchanged across those limits. Operators should therefore evaluate cost limits together with service protections and the resulting restoration decisions. These comparisons can reveal whether additional expenditure enables a fairer plan or whether the protected service requirements restrict the available improvements. The reported outcomes and their operational interpretation apply to the modeled network approximation.

\section{Conclusions}\label{sec:conclusions}
\noindent
This paper develops an optimization framework for incorporating horizontal and vertical fairness into post-disaster power restoration. Using a common operational model with a network approximation, we compare a \utilitarian{} benchmark with the \horizontal{} and \horizontalVertical{} models and examine how different fairness preferences affect repair-crew allocation, power distribution, and restoration cost. The framework also provides complementary measures for evaluating the efficiency and fairness of the resulting restoration plans.

Our results show that fairness in power restoration operates through two interconnected mechanisms: reconnecting isolated areas and allocating limited power within connected areas. The first is largely governed by discrete repair and network decisions and can therefore produce substantial, non-smooth changes in restoration cost. The second provides greater flexibility to redistribute available power among connected loads. This distinction helps explain why stronger fairness preferences do not necessarily translate into proportionally larger efficiency losses. In particular, horizontal fairness can require changes in repair sequencing to improve service for poorly served areas, while vertical fairness can sometimes be improved with limited additional cost by directing scarce power toward more vulnerable communities. These findings highlight the value of considering bus-level service outcomes and community vulnerability alongside system-wide efficiency when planning post-disaster restoration.

Several extensions provide opportunities for future research. First,
the static disaster setting considered here can be extended to dynamic restoration problems in which damage, restoration needs, and resource availability evolve as an event unfolds. Such an extension would also provide a natural setting for incorporating uncertainty in these quantities and updating restoration decisions as new information becomes available. Second, future work could consider additional restoration resources, such as mobile power sources, and their deployment across disconnected areas.

\section*{Data and reproducibility}\label{sec:data}
\ifblindreview
\else
The study data are available on Figshare~\citep{Mao2026figshare} and will be made publicly available after peer review. 
\fi

\ifblindreview
\else
\section*{CRediT authorship contribution statement}
\textbf{Jie Mao:} Data curation, Formal analysis, Investigation, Methodology, Resources, Software, Validation, Visualization, Writing -- original draft, Writing -- review \& editing.

\textbf{Saharnaz Mehrani:} Conceptualization, Methodology, Supervision, Validation, Writing -- review \& editing.

\textbf{Samin Aref:} Project administration, Supervision, Writing -- review \& editing.

\textbf{Payman Dehghanian:} Conceptualization, Writing -- review \& editing.
\fi



\clearpage
\appendix

\section{Proofs of Structural Properties}
\label{app:proofs}

\subsection*{Proof of Lemma~\ref{lem:temporal_invariance}}

By assumption,
\[
\sum_{t\in\mathcal{T}}\gamma_{i,t}
=
\sum_{t\in\mathcal{T}}\tilde{\gamma}_{i,t},
\qquad
\forall i\in\mathcal{B}^L.
\]
It follows from \eqref{eq:yi} that
\[
y_i=\tilde y_i,
\qquad
\forall i\in\mathcal{B}^L.
\]
The \horizontal{} objective depends on a restoration plan only through the vector of average service levels $\{y_i\}_{i\in\mathcal{B}^L}$. Therefore, the two plans have the same \horizontal{} objective value for any $\alpha$. The same argument applies to the MMF objective and, consequently, to DMMF and JI.

For the \horizontalVertical{} objective,
\[
\sum_{t\in\mathcal{T}}(1-\gamma_{i,t})
=
|\mathcal{T}|-\sum_{t\in\mathcal{T}}\gamma_{i,t}.
\]
Hence,
\[
\sum_{t\in\mathcal{T}}(1-\gamma_{i,t})
=
\sum_{t\in\mathcal{T}}(1-\tilde{\gamma}_{i,t}),
\qquad
\forall i\in\mathcal{B}^L.
\]
Multiplying by $SVI_i$ preserves the equality for every bus. The two plans therefore have the same \horizontalVertical{} objective value and VFD. Since $y_i=\tilde y_i$ for every bus, they also have the same VJI. This proves the result.
\hfill\mbox{}\nobreak\hfill$\square$

\subsection*{Proof of Proposition~\ref{prop:connectivity_bound}}

By the definition of $\tau_i$, bus $i$ cannot receive positive service before period $\tau_i$. Therefore,
\[
\gamma_{i,t}=0,
\qquad
t<\tau_i.
\]
Since $0\leq\gamma_{i,t}\leq1$ for all $t\in\mathcal{T}$,
\[
\sum_{t\in\mathcal{T}}\gamma_{i,t}
=
\sum_{t=\tau_i}^{|\mathcal{T}|}\gamma_{i,t}
\leq
|\mathcal{T}|-\tau_i+1.
\]
Using \eqref{eq:yi},
\[
y_i
=
\frac{1}{|\mathcal{T}|}
\sum_{t\in\mathcal{T}}\gamma_{i,t}
\leq
\frac{|\mathcal{T}|-\tau_i+1}{|\mathcal{T}|}.
\]
Thus, for any feasible restoration plan,
\[
\min_{i\in\mathcal{B}^L}y_i
\leq
\min_{i\in\mathcal{B}^L}
\frac{|\mathcal{T}|-\tau_i+1}{|\mathcal{T}|}.
\]
Taking the maximum over $\boldsymbol{\gamma}\in\Gamma_0$ gives the stated upper bound on the MMF service level.

For vertical fairness, $\gamma_{i,t}=0$ for each $t<\tau_i$, which implies
\[
\sum_{t\in\mathcal{T}}(1-\gamma_{i,t})
\geq
\tau_i-1.
\]
Therefore,
\[
SVI_i
\sum_{t\in\mathcal{T}}(1-\gamma_{i,t})
\geq
SVI_i(\tau_i-1),
\qquad
\forall i\in\mathcal{B}^L.
\]
It follows that every feasible restoration plan satisfies
\[
\max_{i\in\mathcal{B}^L}
\left\{
SVI_i
\sum_{t\in\mathcal{T}}(1-\gamma_{i,t})
\right\}
\geq
\max_{i\in\mathcal{B}^L}
\left\{
SVI_i(\tau_i-1)
\right\}.
\]
Taking the minimum over $\boldsymbol{\gamma}\in\Gamma_0$ establishes the second bound.
\hfill\mbox{}\nobreak\hfill$\square$

\subsection*{Proof of Proposition~\ref{prop:fairness_compatibility}}

Suppose a feasible allocation in $\mathcal{Y}^{A}$ provides a minimum service level $m$. The individual service ceilings imply $m\leq\min_i\overline{y}_i$. Moreover, $y_i\geq m$ and $D_i>0$ imply
\[
m\sum_{i\in\mathcal{B}^L}D_i
\leq
\sum_{i\in\mathcal{B}^L}D_i y_i
\leq\overline{E}.
\]
These inequalities give the upper bound in \eqref{eq:aggregate_mmf}. The uniform allocation with every component equal to that bound satisfies all constraints in \eqref{eq:aggregate_service_set}, so the bound is attained.

Next, consider a target vertical-fairness value $z\geq0$. Since $SVI_i>0$, the requirement
\[
\max_{i\in\mathcal{B}^L}
\{SVI_i|\mathcal{T}|(1-y_i)\}
\leq z
\]
is equivalent, together with nonnegativity, to
\[
y_i\geq b_i(z)
:=
\left(1-\frac{z}{|\mathcal{T}|SVI_i}\right)_+,
\qquad \forall i\in\mathcal{B}^L.
\]
The vector $\mathbf{b}(z)$ uses the least energy among allocations meeting this target. It respects the individual ceilings if and only if
\[
z\geq |\mathcal{T}|\max_{i\in\mathcal{B}^L}
\{SVI_i(1-\overline{y}_i)\}.
\]
It respects the energy limit if and only if $\sum_iD_i b_i(z)\leq\overline{E}$. These are precisely the conditions in \eqref{eq:aggregate_vf}. Their feasible set is nonempty, since $z=|\mathcal{T}|\max_iSVI_i$ makes every $b_i(z)$ zero, and closed by continuity. Thus its smallest element exists and equals $V^{A}$.

To attain both benchmarks, an allocation must satisfy $y_i\geq m^{A}$ and $y_i\geq b_i(V^{A})$ for every bus. Since $m^{A}>0$, their joint lower bound is
\[
\ell_i
=
\max\left\{
m^{A},\ 1-\frac{V^{A}}{|\mathcal{T}|SVI_i}
\right\}.
\]
Both constituent lower bounds respect $\overline{y}_i$, so $0\leq\ell_i\leq\overline{y}_i$. Any allocation attaining both benchmarks must use at least $\sum_iD_i\ell_i$ units of energy, which establishes necessity of \eqref{eq:fairness_compatibility}. Conversely, if that condition holds, $\boldsymbol{\ell}\in\mathcal{Y}^{A}$ and it achieves a minimum service level at least $m^{A}$ and a vertical objective at most $V^{A}$. Optimality of the two benchmarks forces both inequalities to hold with equality. This establishes sufficiency.

Finally, suppose $\overline{y}_i=1$ for every bus and $0<\overline{E}<\sum_iD_i$. Then $m^{A}=\overline{E}/\sum_iD_i\in(0,1)$. Any allocation attaining MMF must have $y_i=m^{A}$ for every bus, because increasing any component while keeping all others at least $m^{A}$ would violate the energy limit. If all SVI values equal $w>0$, the vertical objective becomes $w|\mathcal{T}|(1-\min_i y_i)$, so this allocation also attains the vertical-fairness optimum.

If the SVI values are not all equal, write $w_{\max}=\max_iSVI_i$ and let
\[
z_0=|\mathcal{T}|w_{\max}(1-m^{A})>0
\]
be the vertical objective at the uniform MMF allocation. For every bus,
\[
b_i(z_0)
=
\left(1-\frac{w_{\max}}{SVI_i}(1-m^{A})\right)_+
\leq m^{A},
\]
with strict inequality whenever $SVI_i<w_{\max}$. Because every $D_i>0$, it follows that
\[
\sum_iD_i b_i(z_0)
<m^{A}\sum_iD_i
=\overline{E}.
\]
By continuity, there exists $\varepsilon\in(0,z_0)$ such that $\sum_iD_i b_i(z_0-\varepsilon)\leq\overline{E}$. The ceilings equal one, so $\mathbf{b}(z_0-\varepsilon)$ is feasible and achieves a vertical objective at most $z_0-\varepsilon<z_0$. The unique MMF allocation is therefore not vertically optimal, proving that simultaneous attainment is impossible when the SVI values differ.
\hfill\mbox{}\nobreak\hfill$\square$

\subsection*{Proof of Proposition~\ref{prop:restoration_value_structure}}

Fix a combination $k\in\mathcal{K}$ of discrete restoration decisions. Under the power-flow relaxation in Section~\ref{sec:FeasibleRestorationSet}, $\Omega_k$ is convex in the remaining continuous variables. In particular, the line-flow, inverter, and generator capacity constraints are second-order cone constraints, while the voltage-product relaxation is a rotated second-order cone constraint with nonnegative auxiliary voltage variables. The remaining constraints are affine once the discrete decisions are fixed. The service levels $y_i(\mathbf{x})$ and total cost $Z(\mathbf{x})$ are affine, and $V(\mathbf{x})$ is the maximum of affine functions and is therefore convex.

For two feasible cost factors $\lambda_a$ and $\lambda_b$, let $\mathbf{x}^{a}$ and $\mathbf{x}^{b}$ be optimal solutions defining $\Phi_k(\lambda_a)$ and $\Phi_k(\lambda_b)$. For $\eta\in[0,1]$, set
\[
\lambda_{\eta}=\eta\lambda_a+(1-\eta)\lambda_b,
\qquad
\mathbf{x}^{\eta}=\eta\mathbf{x}^{a}+(1-\eta)\mathbf{x}^{b}.
\]
The two plans share the same discrete decisions, so $\mathbf{x}^{\eta}$ retains them and belongs to $\Omega_k$. Because $\beta$ and $\mathbf{y}^{\alpha}$ are fixed,
\[
y_i(\mathbf{x}^{\eta})
=\eta y_i(\mathbf{x}^{a})+(1-\eta)y_i(\mathbf{x}^{b})
\geq\beta y_i^{\alpha},
\qquad \forall i\in\mathcal{B}^L.
\]
Similarly,
\[
Z(\mathbf{x}^{\eta})
=\eta Z(\mathbf{x}^{a})+(1-\eta)Z(\mathbf{x}^{b})
\leq\lambda_{\eta}Z^U.
\]
Thus $\mathbf{x}^{\eta}$ is feasible at $\lambda_{\eta}$. Convexity of $V$ gives
\begin{align*}
\Phi_k(\lambda_{\eta})
&\leq V(\mathbf{x}^{\eta})\\*
&\leq\eta V(\mathbf{x}^{a})+(1-\eta)V(\mathbf{x}^{b})\\*
&=\eta\Phi_k(\lambda_a)+(1-\eta)\Phi_k(\lambda_b).
\end{align*}
This proves convexity. Increasing $\lambda$ only relaxes the cost constraint, so $\Phi_k$ is also nonincreasing. If an infimum is not attained, the same argument applies to arbitrarily near-optimal feasible plans and then takes the approximation error to zero.

Every feasible restoration plan has one of the finitely many discrete combinations in $\mathcal{K}$. Minimizing over these combinations therefore gives $\Phi(\lambda)=\min_{k\in\mathcal{K}}\Phi_k(\lambda)$. It follows that $\Phi$ is nonincreasing as well.

Now suppose that \eqref{eq:schedule_change_condition} holds and that some combination $k$ is optimal at both $\lambda_1$ and $\lambda_3$. Then $\Phi_k(\lambda_1)=\Phi(\lambda_1)$ and $\Phi_k(\lambda_3)=\Phi(\lambda_3)$. Since $\lambda_2=\theta\lambda_1+(1-\theta)\lambda_3$, convexity yields
\begin{align*}
\Phi(\lambda_2)
&\leq\Phi_k(\lambda_2)\\
&\leq\theta\Phi_k(\lambda_1)+(1-\theta)\Phi_k(\lambda_3)\\
&=\theta\Phi(\lambda_1)+(1-\theta)\Phi(\lambda_3),
\end{align*}
contradicting \eqref{eq:schedule_change_condition}. Hence no combination is optimal at both endpoints. In particular, if $\Phi(\lambda_1)=\Phi(\lambda_2)>\Phi(\lambda_3)$, then $0<\theta<1$ makes the right-hand side of \eqref{eq:schedule_change_condition} strictly smaller than $\Phi(\lambda_2)$, establishing the stated implication for a plateau followed by an improvement.
\hfill\mbox{}\nobreak\hfill$\square$

\clearpage
\section{Solver Termination and Objective Bounds}\label{app:solver_certificates}
\setcounter{table}{0}
\renewcommand{\thetable}{B.\arabic{table}}
\renewcommand{\theHtable}{B.\arabic{table}}
\noindent
Tables~\ref{tab:solver_main} and~\ref{tab:solver_identical} report the final solve for each stage, checked against its saved log and solution record. F denotes the primary fairness stage, C the secondary cost stage, and U direct utilitarian cost minimization. HF and MMF maximize their objectives, whereas HVF and cost stages minimize them. O denotes \texttt{OPTIMAL} termination under the solver tolerances. Absolute and relative gaps are $|f-b|$ and $|f-b|/|f|$, computed from the full-precision incumbent $f$ and best bound $b$. Displayed objectives are rounded. Runtimes refer to the individual reported solve.

The main-instance deviation measures use $Z^U=0.1470918846\times10^5$ dollars, $m^{MMF}=0.5833333646$, and $V^{VF}=30$. The identical-demand comparison uses its own utilitarian cost $0.2983568908\times10^5$ dollars. All primary finite-budget HF bounds coincide with their respective incumbent utilities, supporting the plateau comparison within each fixed $\alpha$. For the $\beta=0$ application of Proposition~\ref{prop:restoration_value_structure}, the middle lower bound is 60 and feasible endpoint objectives are 65 and 45. Thus, $L_2-0.4U_1-0.6U_3=7$, well beyond numerical tolerances.

\begingroup
\singlespacing
\setlength{\tabcolsep}{3pt}
\fontsize{10}{12}\selectfont
\begin{longtable}{@{}llrrrrrc@{}}
\caption{Main-instance solver results. Cost objectives are in units of $10^5$ dollars.}\label{tab:solver_main}\\
\toprule
Case, $\lambda$ & Stage & Incumbent & Bound & Abs. gap & Rel. gap & \shortstack{Time\\(s)} & Status\\
\midrule
\endfirsthead
\multicolumn{8}{l}{\tablename\ \thetable\ continued}\\
\toprule
Case, $\lambda$ & Stage & Incumbent & Bound & Abs. gap & Rel. gap & \shortstack{Time\\(s)} & Status\\
\midrule
\endhead
\midrule\multicolumn{8}{r}{Continued on next page}\\\endfoot
\bottomrule\endlastfoot
UT, $\infty$ & U & 0.147091885 & 0.147091885 & 0 & 0 & 9.29 & O \\
HF, $\alpha=0.25$, $\infty$ & F & 39.4203647 & 39.4203647 & 4.38e-10 & 1.11e-11 & 11.54 & O \\
HF, $\alpha=0.25$, $\infty$ & C & 0.22657499 & 0.226574983 & 6.94e-09 & 3.06e-08 & 4.68 & O \\
HF, $\alpha=0.5$, $\infty$ & F & 60.5444924 & 60.5444924 & 0 & 0 & 8.04 & O \\
HF, $\alpha=0.5$, $\infty$ & C & 0.226574975 & 0.226574975 & 3.58e-11 & 1.58e-10 & 2.65 & O \\
HF, $\alpha=0.75$, $\infty$ & F & 124.306504 & 124.306504 & 0 & 0 & 8.88 & O \\
HF, $\alpha=0.75$, $\infty$ & C & 0.226572336 & 0.226572336 & 0 & 0 & 4.13 & O \\
HF, $\alpha=1$, $\infty$ & F & -3.9662085 & -3.9662085 & 0 & 0 & 5.43 & O \\
HF, $\alpha=1$, $\infty$ & C & 0.226572377 & 0.226572377 & 0 & 0 & 5.00 & O \\
HF, $\alpha=2$, $\infty$ & F & -37.3213277 & -37.3213277 & 0 & 0 & 30.89 & O \\
HF, $\alpha=2$, $\infty$ & C & 0.276101935 & 0.276101028 & 9.07e-07 & 3.29e-06 & 28.60 & O \\
HF, $\alpha=5$, $\infty$ & F & -18.0948734 & -18.0948734 & 0 & 0 & 6.28 & O \\
HF, $\alpha=5$, $\infty$ & C & 0.257411212 & 0.257411212 & 0 & 0 & 5.06 & O \\
HF, $\alpha=10$, $\infty$ & F & -45.7798324 & -45.7798324 & 0 & 0 & 2.89 & O \\
HF, $\alpha=10$, $\infty$ & C & 0.257411212 & 0.257411212 & 4.64e-11 & 1.80e-10 & 9.75 & O \\
MMF, $\infty$ & F & 0.583333365 & 0.583333365 & 0 & 0 & 4.19 & O \\
MMF, $\infty$ & C & 0.224101114 & 0.224101114 & 0 & 0 & 4.30 & O \\
HF, $\alpha=1$, 1.25 & F & -4.85057039 & -4.85057039 & 0 & 0 & 2.29 & O \\
HF, $\alpha=1$, 1.25 & C & 0.159795984 & 0.159795984 & 1.11e-14 & 6.95e-14 & 2.58 & O \\
HF, $\alpha=2$, 1.25 & F & -39.075302 & -39.075302 & 0 & 0 & 1.81 & O \\
HF, $\alpha=2$, 1.25 & C & 0.159795984 & 0.159795984 & 0 & 0 & 2.93 & O \\
HF, $\alpha=10$, 1.25 & F & -1181.17243 & -1181.17243 & 0 & 0 & 2.25 & O \\
HF, $\alpha=10$, 1.25 & C & 0.159795984 & 0.159795984 & 0 & 0 & 2.68 & O \\
HF, $\alpha=1$, 1.4 & F & -4.85057039 & -4.85057039 & 0 & 0 & 3.84 & O \\
HF, $\alpha=1$, 1.4 & C & 0.159795984 & 0.159795984 & 0 & 0 & 1.84 & O \\
HF, $\alpha=2$, 1.4 & F & -39.075302 & -39.075302 & 0 & 0 & 9.66 & O \\
HF, $\alpha=2$, 1.4 & C & 0.159795984 & 0.159795984 & 0 & 0 & 3.86 & O \\
HF, $\alpha=10$, 1.4 & F & -1181.17243 & -1181.17243 & 0 & 0 & 2.86 & O \\
HF, $\alpha=10$, 1.4 & C & 0.159795984 & 0.159795984 & 0 & 0 & 2.88 & O \\
HF, $\alpha=1$, 1.5 & F & -4.85057039 & -4.85057039 & 0 & 0 & 5.79 & O \\
HF, $\alpha=1$, 1.5 & C & 0.159795984 & 0.159795984 & 0 & 0 & 2.75 & O \\
HF, $\alpha=2$, 1.5 & F & -39.075302 & -39.075302 & 0 & 0 & 13.70 & O \\
HF, $\alpha=2$, 1.5 & C & 0.159795984 & 0.159795984 & 0 & 0 & 5.67 & O \\
HF, $\alpha=10$, 1.5 & F & -1181.17243 & -1181.17243 & 0 & 0 & 5.46 & O \\
HF, $\alpha=10$, 1.5 & C & 0.159795984 & 0.159795984 & 0 & 0 & 2.21 & O \\
HF, $\alpha=0.1$, $\infty$ & F & 32.422741 & 32.422741 & 9.00e-09 & 2.78e-10 & 6.46 & O \\
HF, $\alpha=0.1$, $\infty$ & C & 0.226574858 & 0.226574858 & 7.97e-11 & 3.52e-10 & 5.35 & O \\
HVF, $\beta=0$, $\infty$ & F & 30 & 30 & 0 & 0 & 2.88 & O \\
HVF, $\beta=0$, $\infty$ & C & 0.25940236 & 0.25940236 & 0 & 0 & 6.51 & O \\
HVF, $\beta=0.75$, $\infty$ & F & 30 & 30 & 0 & 0 & 2.89 & O \\
HVF, $\beta=0.75$, $\infty$ & C & 0.259668764 & 0.259668764 & 0 & 0 & 3.09 & O \\
HVF, $\beta=0.85$, $\infty$ & F & 30 & 30 & 0 & 0 & 2.11 & O \\
HVF, $\beta=0.85$, $\infty$ & C & 0.259668767 & 0.259668765 & 1.45e-09 & 5.57e-09 & 3.23 & O \\
HVF, $\beta=0.9$, $\infty$ & F & 30 & 30 & 0 & 0 & 2.43 & O \\
HVF, $\beta=0.9$, $\infty$ & C & 0.273487462 & 0.273487027 & 4.35e-07 & 1.59e-06 & 3.82 & O \\
HVF, $\beta=0$, 1.25 & F & 64.9999999 & 64.9999999 & 0 & 0 & 3.26 & O \\
HVF, $\beta=0$, 1.25 & C & 0.155307549 & 0.155307534 & 1.59e-08 & 1.03e-07 & 11.75 & O \\
HVF, $\beta=0.75$, 1.25 & F & 70 & 70 & 0 & 0 & 2.31 & O \\
HVF, $\beta=0.75$, 1.25 & C & 0.159795984 & 0.159795984 & 0 & 0 & 1.91 & O \\
HVF, $\beta=0$, 1.4 & F & 60 & 60 & 0 & 0 & 3.45 & O \\
HVF, $\beta=0$, 1.4 & C & 0.191796899 & 0.191796882 & 1.65e-08 & 8.58e-08 & 10.59 & O \\
HVF, $\beta=0.75$, 1.4 & F & 70 & 70 & 0 & 0 & 2.25 & O \\
HVF, $\beta=0.75$, 1.4 & C & 0.159795984 & 0.159795984 & 0 & 0 & 1.92 & O \\
HVF, $\beta=0$, 1.5 & F & 45 & 45 & 0 & 0 & 19.26 & O \\
HVF, $\beta=0$, 1.5 & C & 0.219800311 & 0.219800311 & 0 & 0 & 4.28 & O \\
HVF, $\beta=0.75$, 1.5 & F & 70 & 70 & 0 & 0 & 2.08 & O \\
HVF, $\beta=0.75$, 1.5 & C & 0.159795984 & 0.159795984 & 1.09e-14 & 6.81e-14 & 1.54 & O \\

\end{longtable}
\endgroup

\begingroup
\singlespacing
\setlength{\tabcolsep}{3pt}
\fontsize{10}{12}\selectfont
\begin{longtable}{@{}llrrrrrc@{}}
\caption{Identical-demand solver results. Cost objectives are in units of $10^5$ dollars.}\label{tab:solver_identical}\\
\toprule
Case, $\lambda$ & Stage & Incumbent & Bound & Abs. gap & Rel. gap & \shortstack{Time\\(s)} & Status\\
\midrule
\endfirsthead
\multicolumn{8}{l}{\tablename\ \thetable\ continued}\\
\toprule
Case, $\lambda$ & Stage & Incumbent & Bound & Abs. gap & Rel. gap & \shortstack{Time\\(s)} & Status\\
\midrule
\endhead
\midrule\multicolumn{8}{r}{Continued on next page}\\\endfoot
\bottomrule\endlastfoot
UT, $\infty$ & U & 0.298356891 & 0.298356888 & 3.05e-09 & 1.02e-08 & 11.76 & O \\
MMF, $\infty$ & F & 0.583333333 & 0.583333333 & 0 & 0 & 4.20 & O \\
MMF, $\infty$ & C & 0.466645747 & 0.466645747 & 0 & 0 & 1.45 & O \\
HF, $\alpha=0.1$, $\infty$ & F & 30.2053221 & 30.2053221 & 6.44e-08 & 2.13e-09 & 79.59 & O \\
HF, $\alpha=0.1$, $\infty$ & C & 0.298356524 & 0.298356524 & 0 & 0 & 3.52 & O \\
HF, $\alpha=0.25$, $\infty$ & F & 37.0525014 & 37.0525019 & 5.64e-07 & 1.52e-08 & 3.98 & O \\
HF, $\alpha=0.25$, $\infty$ & C & 0.298356812 & 0.298356812 & 0 & 0 & 2.66 & O \\
HF, $\alpha=0.5$, $\infty$ & F & 57.9265633 & 57.9265633 & 0 & 0 & 3.31 & O \\
HF, $\alpha=0.5$, $\infty$ & C & 0.301007008 & 0.301007008 & 0 & 0 & 2.30 & O \\
HF, $\alpha=0.75$, $\infty$ & F & 121.432778 & 121.432778 & 0 & 0 & 8.92 & O \\
HF, $\alpha=0.75$, $\infty$ & C & 0.301007044 & 0.301007008 & 3.59e-08 & 1.19e-07 & 3.24 & O \\
HF, $\alpha=1$, $\infty$ & F & -7.14300167 & -7.14300167 & 0 & 0 & 3.88 & O \\
HF, $\alpha=1$, $\infty$ & C & 0.301007041 & 0.301007041 & 0 & 0 & 3.68 & O \\
HF, $\alpha=2$, $\infty$ & F & -41.8354181 & -41.8354172 & 9.09e-07 & 2.17e-08 & 15.92 & O \\
HF, $\alpha=2$, $\infty$ & C & 0.317298081 & 0.317298081 & 0 & 0 & 2.02 & O \\
HF, $\alpha=5$, $\infty$ & F & -39.2471887 & -39.2471887 & 0 & 0 & 20.27 & O \\
HF, $\alpha=5$, $\infty$ & C & 0.466645747 & 0.466645747 & 0 & 0 & 3.93 & O \\

\end{longtable}
\endgroup

\end{document}